\documentclass[smallextended,envcountsect]{svjour3}
\smartqed  

\usepackage{mathtools}
\usepackage{graphicx}
\usepackage{amssymb}
\usepackage{url}
\usepackage{multirow}
\usepackage{algorithm}
\usepackage{algpseudocode}
\usepackage{booktabs}
\usepackage[dvipsnames]{xcolor}
\usepackage{colortbl}
\usepackage{tcolorbox}
\tcbuselibrary{breakable}
\usepackage{hyperref}       

\colorlet{mypurple}{purple!20}
\usepackage[color=mypurple]{todonotes}

\newenvironment{thm}{\begin{theorem}}{\end{theorem}}

\newenvironment{cor}{\begin{corollary}}{\end{corollary}}
\newenvironment{prop}{\begin{proposition}}{\end{proposition}}
\newenvironment{bem}{\begin{remark}}{\end{remark}}

\spnewtheorem{obs}{Observation}[section]{\bfseries}{\itshape}

\newcounter{rkcounter}
\renewcommand{\therkcounter}{RK.\arabic{rkcounter}}

\newcommand{\RK}[2][]{%
  \refstepcounter{rkcounter}%
  \begin{tcolorbox}[
    colback=Peach, 
    colframe=black, 
    size=small, 
    title=\textbf{\therkcounter}\ifx&#1&\else\ - #1\fi
  ]%
    
    \small #2
  \end{tcolorbox}%
}

\newcommand{\cardy}[1]{{\mathbf{e}}^\top{#1}}
\newcommand{\bqp}[2]{\frac{1}{2}{#2}^\top{#1}#2}

\newcommand{\inrnn}[1]{#1 \in \mathbb{R}^{n \times n}}
\newcommand{\inbin}[1]{#1 \in \{0, 1\}^n}
\newcommand{\linpen}[2]{\lambda(#1-#2)}
\newcommand{\qpen}[2]{\eta(#1-#2)^2}
\newcommand{\aval}[1]{ A_k := \max_{0 \leq l<k}\left\{ \frac{#1_{k}-#1_{l}}{k-l}\right\}}
\newcommand{\bval}[1]{ B_k := \min_{k<l \leq n}\left\{ \frac{#1_{l}-#1_{k}}{l-k}\right\}}
\newcommand{\valapprox}[1]{\tilde{#1}_k}

\newcommand{\lambdaapprox}{\frac{\tilde{A}_k+\tilde{B}_k}{2}}
\newcommand{\etaapprox}{\frac{\tilde{A}_k-\tilde{B}_k}{2}}

\makeatletter
\renewcommand{\makeheadbox}{}
\makeatother
\renewcommand{\subclassname}{{\bfseries Mathematics Subject Classification}\enspace}

\title{Linear and Quadratic Penalties for Cardinality Constraints in
  Unconstrained Binary Optimization}
\titlerunning{Linear and Quadratic Penalties for Cardinality Constraints}
\authorrunning{Ku\v{z}el et al.}
\author{Roman~Ku\v{z}el \and
        M\'{a}ty\'{a}s~Koniorczyk \and
        Janez~Povh}
\institute{R.~Ku\v{z}el \at
              Rudolfovo -- Science and Technology Centre Novo Mesto,
              Novo Mesto, Slovenia \\
              \email{roman.kuzel@rudolfovo.eu}
      \and
           M.~Koniorczyk, Corresponding author \at
              HUN-REN Wigner Research Centre for Physics,
              Budapest, Hungary \\
              \email{koniorczyk.matyas@wigner.hun-ren.hu}
      \and
           J.~Povh \at
              Rudolfovo -- Science and Technology Centre Novo Mesto,
              Novo Mesto, Slovenia \\
              \email{janez.povh@rudolfovo.eu}
             }
\date{20 November 2026}

\begin{document}

\maketitle

\begin{abstract}
  We study penalty strategies for converting constrained binary
  optimization problems to quadratic unconstrained binary optimization
  (QUBO) form on the example of the sparsest $k$-subgraph problem. We
  adopt an augmented Lagrangian approach involving a linear penalty in
  addition to the commonly used quadratic one, and show the benefits
  of this approach. From the constraint structure, we derive
  sufficient bounds on the penalty coefficients. We conduct extensive
  experiments comparing combinations of our different penalty
  strategies with a wide range of classical and quantum QUBO
  solvers. Our experiments show that to obtain a good-quality
  solution, the right penalty strategy is almost as important as the
  quality of the QUBO solver, and iterative approaches appear to be
  the most efficient. We discuss the mechanisms that lead to a
  degradation of the solution quality due to large quadratic
  penalties.  

\end{abstract}

\section{Introduction}
\label{sec:introduction}

Quadratic unconstrained binary optimization (QUBO)~\cite{Punnen_2022}
has attracted renewed research interest because of its relevance in
quantum computing.  Being equivalent to the Ising problem, QUBO is the
native input format of quantum annealers~\cite{Kadowaki_1998,Das_2008}
such as D-Wave~\cite{Johnson_2011}, and it can also serve as the input
of algorithms for gate-based quantum computers, such as
QAOA~\cite{Farhi_2014,Blekos_2024}.  QUBO is also equivalent to
Max-Cut, for which many good classical heuristics are
available~\cite{MQLIB}.  Exact solvers such as BiqBin~\cite{BIQBIN}
and Gurobi~\cite{Gurobi} have developed significantly, and new
physics-based algorithms such as simulated
bifurcation~\cite{Goto_2019,SBM,Kanao_2022,2604.01050} offer
quantum-like behavior in some cases.

Many real-world applications, however, lead to constrained problems.
To use the above solvers, such a problem has to be transformed to QUBO
form, most commonly by quadratic penalties: the square of the
violation of an equality constraint, multiplied by a large enough
number, is added to the objective~\cite{Lucas_2014,Glover_2022}.  The
appropriate choice of the large enough multiplier is often not
trivial.

Lasserre~\cite{Lasserre2016} and Gusmeroli and Wiegele~\cite{EXPEDIS}
derived general penalty bounds for linearly constrained quadratic
binary programs with integer coefficients. The penalties obtained in this way
guarantee optimality, but they can be hard to calculate and are often
too big in practice. Hence, it is worth taking into account the
structure of the constraints in the actual problem. For the one-hot
encodings of permutation constraints, for instance,
Ayodele~\cite{Ayodele_2022} derives closed-form static penalties from
single-flip changes of the QUBO, which outperform global bounds on
Fujitsu's Digital Annealer. These bounds, however, are not provably
tight.

Several empirical approaches tune the penalties for specific problem
classes in the quantum setting, by the cross-entropy
method~\cite{Roch_2021}, by machine-learning
regression~\cite{Roch_2023}, or by characterizing the penalties of
multi-knapsack variants under QAOA and quantum
annealing~\cite{Guney_2025}.  These and most other applied
contributions in the quantum realm are specific to quantum solvers,
and typically address small problems because of the limited size of
the available hardware.  Quintero et al.~\cite{Quintero_2022} study
QUBO formulations of the maximum $k$-colorable subgraph problem in
detail.  Building on the known formulations of the max-clique and
stable set problems, they introduce linear-based and nonlinear QUBO
reformulations and, exploiting the properties of the constraints,
analyze the optimal choice of the penalty parameters. They show the
benefits of a well-chosen encoding and penalty values through the time
dependence of the energy gap during annealing, the embedding
properties, and the time-to-solution on different D-Wave processors.
Vysko\v{c}il and Djidjev~\cite{Vyskocil_2019} consider a cardinality
constraint whose squared penalty makes the problem graph complete, so
that its minor embedding requires a number of physical qubits
quadratic in~$n$.  They replace it by a sparse Ising program assembled
from small unit-cell gadgets, maximizing the gap between feasible and
infeasible assignments under the coefficient range of the hardware.
How the objective and the constraint share this range is left open.

Large penalties can also be avoided by Lagrangian methods.  Ronagh,
Woods, and Iranmanesh~\cite{10.5555/3179330.3179336} use a
branch-and-bound algorithm with Lagrangian lower bounds, and Karimi and
Ronagh~\cite{Karimi_2017} a subgradient method for the Lagrangian dual
of the stable set problem, demonstrated on the D-Wave 2X.
Djidjev~\cite{Djidjev_2023a} applies the augmented Lagrangian method
directly to the QUBO encoding of set cover, as we do, with iteratively
updated linear and quadratic multipliers, but the updates are purely
heuristic, with no bound on the required penalty or on the resulting
distortion; in a companion work he applies the method to the chain
couplings of a minor embedding~\cite{Djidjev_2023b}.  Cellini et
al.~\cite{Cellini_2024} derive a closed-form estimate of the
multipliers for bin packing, similar in spirit to our sufficiency
bounds but heuristic rather than provably tight, and Liu and
Moraglio~\cite{Liu_2025} search for the penalty coefficients
iteratively in a multi-constraint setting with an empirical sensitivity
heuristic.

The general bounds thus do not exploit the structure of the
constraints, the augmented Lagrangian approaches estimate the
multipliers heuristically, and the empirical studies are mostly
specific to quantum solvers. In the present contribution we
demonstrate that exploiting the problem structure is crucial for
finding the right encoding, and present considerations that contribute
to the understanding of the drawbacks of large quadratic penalties.
We showcase this on the sparsest $k$-subgraph problem, the complement
of the well-studied densest $k$-subgraph
problem~\cite{Feige_2001,Hanaka_2023}, which is relevant in many
fields including random forests~\cite{Kisander_prep}.  It has a single
equality constraint, and the nature of its violations allows for
adding a linear penalty term and adopting the augmented Lagrangian
method.  Using a broad range of classical, quantum, and
quantum-inspired heuristics, we demonstrate
the benefits of linear penalties, and show empirically that iterative
methods can be very useful, especially in combination with certain
solvers.
This allows us to study how the penalty strategy and the choice of
solver interact in determining solution quality and computational
cost, in particular to what extent a well-chosen penalty strategy can
compensate for a weaker solver, and vice versa.

Section~\ref{sec:card_qbp} derives conditions under which the QUBO
encoding of a cardinality-constrained quadratic binary problem is
exact, Section~\ref{sec:penalties_solvers} discusses how large
quadratic penalties affect inner solvers, and
Section~\ref{sec:application_sks} applies these results to the sparsest
$k$-subgraph problem.  Sections~\ref{sec:computational_study}
and~\ref{sec:computational_results} present the computational study.
Section~\ref{sec:cs_cor} compares our penalty bounds with
those in the literature, and Section~\ref{sec:cs_conclusion} summarizes
our results.

\section{Quadratic Binary Problems with a Cardinality Constraint}
\label{sec:card_qbp}

Let us consider a quadratic binary optimization problem with a constraint of a
specific form:
\begin{align}\label{rk:QBO_constrainedl}
\min \left \{\bqp{Q}{x}\mid \inbin{x}, \cardy{x}=k, \inrnn{Q} \right \},
\end{align}
that is, the cardinality (or Hamming weight) of the solution vector
should be equal to $k$.  We assume $1 \leq k \leq n-1$; for
$k \in \{0, n\}$ the feasible set consists of a single point and the
problem is trivial.  We consider $Q$ to be symmetric which can be
done wlog., and the factor $1/2$ is for notational convenience.  In order to use a solver that is capable of solving QUBO problems (which will
be termed as \emph{inner solver} throughout this paper) to solve this
constrained problem, one has to encode it into a QUBO
\begin{align}\label{rk:QUBO}
\min \left \{f_{\lambda,\eta}(x) \mid \inbin{x} \right \},
\end{align}
where
\begin{equation*}
f_{\lambda,\eta}(x) = \bqp{Q}{x} -\linpen{\cardy{x}}{k}+\qpen{\cardy{x}}{k},
\end{equation*}
and $\mathbf{e}$ is the all-1 vector.
In this section we derive conditions under which
problem~\eqref{rk:QBO_constrainedl} is equivalent to the
QUBO~\eqref{rk:QUBO}.

The objective of the QUBO contains a quadratic penalty term $\eta$ and a
Lagrangian multiplier $\lambda$; our conditions will apply to these terms.
Our motivation is that this reformulation enables one to keep the magnitude of
the quadratic penalty term under control, which is essential to maintain the
stability of the inner solver.

We use $*$ to indicate a minimizer of an optimization problem, i.e.,
$x^*$ is a binary vector minimizing the problem under discussion.
Minimizers wit cardinality $l$ will be denoted by $x_l^*$, and $x^*$ will be understood as $x_k^*$.
For $l \in \{0, 1, \dots, n\}$, let $q_l$ be the minimum of the constrained
problem~\eqref{rk:QBO_constrainedl}, setting the cardinality in
the constraint to $l$. Now let us introduce the following quantities:
\begin{align}\label{rk:ABvals}
\aval{q},\quad
\bval{q}
\end{align}

We can now proceed to the main result of this work showing tight dependence of
both linear penalty term $\lambda$ and quadratic penalty term $\eta$ on these
quantities when transforming the constrained binary optimization problem
\eqref{rk:QBO_constrainedl} to its equivalent unconstrained version
\eqref{rk:QUBO}.

\begin{thm}\label{rk:theorem}
Let $Q \in \mathbb{R}^{n \times n}$ and let $\valapprox{A}, \valapprox{B}$ be such that
$\valapprox{A}\geq A_k, \valapprox{B}\leq B_k$, and let $\epsilon$ be an arbitrary positive number. Then the following statements are equivalent:
\begin{enumerate}
\item[(a)]\label{rk:one} $x^*$ is a minimizer of \eqref{rk:QBO_constrainedl}.
\item[(b)]\label{rk:two} $x^*$ is a minimizer of \eqref{rk:QUBO}
  for any $(\lambda, \eta)$ satisfying
  \begin{align}
    &\eta \geq\max\left (0, \etaapprox+\epsilon\right ),\quad \epsilon>0, \nonumber \\
    &\lambda\in \left[\lambdaapprox-\epsilon',\quad \lambdaapprox+\epsilon'\right],\
    \epsilon' \geq 0, \epsilon'<\epsilon. \label{eq:lambdaetacond}
  \end{align}
\end{enumerate}
\end{thm}
\begin{proof}
$(a)\Rightarrow (b)$:\\
Let $x_k^*$ be a minimizer of \eqref{rk:QBO_constrainedl} for some $k$. By
contradiction, assume that $f_{\lambda,\nu}(x^*_k)$ is not equal to the optimum of \eqref{rk:QUBO} for $\lambda$,
$\eta$ and $l$, i.e.:
\[
q_k\ge q_{l} - \linpen{l}{k} + \qpen{l}{k}.
\]
Without loss of generality, assume that $l<k$. We can rewrite the above
inequality as follows:
\[
\frac{q_k-q_{l}}{k-l}\ge\lambda+\eta(k-l)
\]
Since $l<k$, definition \eqref{rk:ABvals} implies that
\[
A_k\geq\lambda+\eta.
\]
\[
A_k\geq\lambdaapprox-\epsilon'+\etaapprox+\epsilon.
\]
Since $\epsilon'<\epsilon$ we have
\[A_k>\valapprox{A},\] a contradiction.\\
$(a) \Leftarrow (b)$:\\
Assume that $x_{l}^*:=x^*$ is a minimizer of \eqref{rk:QUBO} such that
$\cardy{x^*}=l$.
Then
\[
q_{l}-\linpen{l}{k}+\qpen{l}{k}\leq q_{l'}-\linpen{l'}{k}+\qpen{l'}{k},
\quad \forall l'\in \{0, 1,\dots,n\}.
\]
For $l'=k$ we have
\[
q_{l}-\linpen{l}{k}+\qpen{l}{k}\leq q_{k}
\]
Without loss of generality let $l<k$. Then we have
\[
\lambda+\eta(k-l)\leq\frac{q_k-q_{l}}{k-l}.
\]
Using the same arguments as above we get
\[
\valapprox{A}<A_k.
\]
This contradicts assumptions of (b). Therefore for any minimizer $x^*$ of
\eqref{rk:QUBO} we have that $\cardy{x^*}=k$ and $x^*$ is also minimizer of
\eqref{rk:QBO_constrainedl}. \qed
\end{proof}

A straightforward corollary of the above Theorem is the following.
\begin{cor}
For any $(\lambda,\eta)$ satisfying~\eqref{eq:lambdaetacond}, any
minimizer $x^*$ of \eqref{rk:QUBO} has cardinality $k$.
\end{cor}
A further consequence of Theorem~\ref{rk:theorem} is that when $B_k > A_k$ the
quadratic penalty can be dropped entirely: the cardinality constraint is
enforced by the linear term alone:
\begin{cor}\label{rk:lambda}
If $B_k-A_k>0$ for values from \eqref{rk:ABvals}, then the constrained
quadratic binary problem \eqref{rk:QBO_constrainedl} is equivalent to
unconstrained quadratic binary problem \eqref{rk:QUBO} with the linear penalty
term $\lambda\in(A_k, B_k)$ and the quadratic penalty term $\eta=0$.
\end{cor}
\begin{proof}
We use $\valapprox{A}:=A_k$ and $\valapprox{B}:=B_k$. Since $B_k-A_k > 0$ we
can choose $0\leq\epsilon'<\frac{B_k-A_k}{2}$. Therefore,
$\frac{A_k-B_k}{2}+\epsilon' < 0\implies\eta=0$. \qed
\end{proof}

At first sight a linear penalty $-\lambda(\mathbf{e}^\top x - k)$ seems
inadequate for enforcing an \emph{equality} constraint: unlike the quadratic
term $\eta(\mathbf{e}^\top x-k)^2\geq 0$, which is non-negative and
symmetrically penalizes violations in both directions, the linear term is
signed and for $\lambda>0$ rewards solutions with $\mathbf{e}^\top x > k$.
The corollary shows that the right $\lambda$ nevertheless enforces the
constraint from both sides via the structure of the optimal value function
$q_l$: choosing $\lambda\in(A_k,B_k)$ simultaneously ensures
$q_l - \lambda(l-k) > q_k$ for all $l > k$ (adding elements is too costly)
and for all $l < k$ (removing elements is too costly), so the unconstrained
minimum is trapped at cardinality $k$.  Positive semidefiniteness of the
penalty is not required; what is required is the non-empty interval
$(A_k,B_k)$, i.e., $B_k > A_k$.

Corollary~\ref{rk:lambda} gives a sufficient condition for the linear
penalty alone to enforce the constraint.  The following proposition
strengthens this to a complete characterization: the quadratic penalty can
be set to zero if and only if the interval $(A_k,B_k)$ is non-empty.

\begin{prop}\label{prop:eta_zero_iff}
The QUBO~\eqref{rk:QUBO} with $\eta=0$ is equivalent
to~\eqref{rk:QBO_constrainedl} if and only if\/ $A_k \leq B_k$ and
$\lambda \in (A_k,\,B_k)$.
\end{prop}
\begin{proof}
For~\eqref{rk:QUBO} with $\eta=0$ to be equivalent
to~\eqref{rk:QBO_constrainedl}, we need
\begin{enumerate}
  \item $q_{k_-} + \lambda(k-k_-) \geq q_k$ for all $k_-<k$,
  \item $q_{k_+} - \lambda(k_+-k) \geq q_k$ for all $k_+>k$.
\end{enumerate}
These are equivalent to $\lambda \geq A_k$ and $\lambda \leq B_k$,
respectively.  Both can be satisfied simultaneously if and only if
$A_k \leq B_k$, with admissible values $\lambda \in (A_k,\,B_k)$. \qed
\end{proof}

When $A_k > B_k$ a positive quadratic penalty is indispensable, and a
natural question is how small it can be made.  To state the answer we
introduce the \emph{unit-step} analogues of $A_k$ and $B_k$, which
restrict the optimization in~\eqref{rk:ABvals} to the immediately
neighboring cardinalities:

\begin{definition}\label{def:unit_step}
For $1 \leq k \leq n-1$, let
\begin{enumerate}
  \item $A_k^1 := q_k - q_{k-1}$,
  \item $B_k^1 := q_{k+1} - q_k$.
\end{enumerate}
\end{definition}

Note that by definition $A_k \geq A_k^1$ and $B_k \leq B_k^1$: the
global quantities $A_k,B_k$ are more demanding than their unit-step
counterparts.

\begin{thm}\label{thm:min_eta}
Let $\eta^*$ be the infimum of values of\/ $\eta\ge0$ for which
\eqref{rk:QBO_constrainedl} and~\eqref{rk:QUBO} are equivalent.
Then
\begin{enumerate}
  \item[(a)] if\/ $A_k \leq B_k$: \;$\eta^*=0$;\; equivalence
    holds for $\eta = 0$ and any $\lambda \in (A_k,\,B_k)$
    when $A_k < B_k$, and for any $\eta > 0$ with
    $\lambda = A_k$ when $A_k = B_k$;
  \item[(b)] if\/ $A_k > B_k$:
    \;$\eta^* \in \bigl[\tfrac{A_k^1-B_k^1}{2},\;
    \tfrac{A_k-B_k}{2}\bigr]$ and equivalence holds for any
    $\eta > \eta^*$.
\end{enumerate}
\end{thm}

\begin{proof}
\emph{(a)} follows immediately from
Proposition~\ref{prop:eta_zero_iff}.

\emph{(b)}
The upper bound $\eta^* \leq (A_k-B_k)/2$ follows from
Theorem~\ref{rk:theorem}.  For the lower bound, consider the
necessary conditions at cardinalities $k\pm 1$:
\begin{align*}
  q_{k-1} + \lambda + \eta &\geq q_k,
    \\
  q_{k+1} - \lambda + \eta &\geq q_k.
\end{align*}
Using Definition~\ref{def:unit_step}, these become
$\eta \geq A_k^1 - \lambda$ and $\lambda \leq \eta + B_k^1$.  To
minimize $\eta$ we maximize $\lambda$; the second inequality then gives
$\lambda \leq \eta + B_k^1$.  Substituting into the first:
$\eta \geq A_k^1 - \eta - B_k^1$, whence
$\eta \geq (A_k^1-B_k^1)/2$.  The bounds on $\lambda$ follow by the
same chain of inequalities. \qed
\end{proof}

\begin{cor}\label{cor:eta_monotone}
For any $\eta > \eta^*$, every minimizer of~\eqref{rk:QUBO} has
cardinality~$k$.
\end{cor}

The conditions~\eqref{eq:lambdaetacond} of Theorem~\ref{rk:theorem}
imply the strict inequalities $\lambda+\eta > A_k$ and
$\lambda-\eta < B_k$; it is the parameter $\epsilon > 0$ that makes
them strict.  The following proposition describes what remains true at
the boundary, that is, for $\epsilon = 0$.

\begin{prop}[Boundary case]\label{prop:boundary}
Let $\eta > 0$, $\lambda+\eta \geq A_k$ and $\lambda-\eta \leq B_k$.
Then
\begin{enumerate}
  \item[(a)] the optimal value of~\eqref{rk:QUBO} is $q_k$, and every
    minimizer of~\eqref{rk:QBO_constrainedl} is a minimizer
    of~\eqref{rk:QUBO};
  \item[(b)] every minimizer of~\eqref{rk:QUBO} of cardinality~$k$ is a
    minimizer of~\eqref{rk:QBO_constrainedl};
  \item[(c)] every other minimizer of~\eqref{rk:QUBO} has cardinality
    $k-1$ or $k+1$.  Minimizers of cardinality $k-1$ exist if and only
    if $A_k^1 = \lambda+\eta$, and minimizers of cardinality $k+1$ exist
    if and only if $B_k^1 = \lambda-\eta$.
\end{enumerate}
\end{prop}
\begin{proof}
For $l \in \{0,1,\dots,n\}$ let
$\varphi(l) := q_l - \linpen{l}{k} + \qpen{l}{k}$ be the minimum of
$f_{\lambda,\eta}$ over the vectors of cardinality~$l$; in particular
$\varphi(k) = q_k$.  For $l<k$ and $d := k-l$,
\[
  \varphi(l) - q_k
  \;=\; d\Bigl(\lambda + \eta d - \frac{q_k-q_l}{d}\Bigr)
  \;\geq\; d\bigl(\lambda+\eta-A_k + \eta(d-1)\bigr) \;\geq\; 0,
\]
with equality if and only if $d=1$ and $q_k-q_{k-1} = \lambda+\eta$.
For $l>k$ and $d := l-k$,
\[
  \varphi(l) - q_k
  \;=\; d\Bigl(\frac{q_l-q_k}{d} - \lambda + \eta d\Bigr)
  \;\geq\; d\bigl(B_k-(\lambda-\eta) + \eta(d-1)\bigr) \;\geq\; 0,
\]
with equality if and only if $d=1$ and $q_{k+1}-q_k = \lambda-\eta$.
Hence $\min_l \varphi(l) = q_k$, which gives (a) and (c).  For (b), a
minimizer $x$ of~\eqref{rk:QUBO} with $\cardy{x}=k$ satisfies
$\bqp{Q}{x} = f_{\lambda,\eta}(x) = q_k$. \qed
\end{proof}

Claim (b) does not depend on the assumptions on $\lambda$ and $\eta$:
on vectors of cardinality~$k$ the function $f_{\lambda,\eta}$ coincides
with the objective of~\eqref{rk:QBO_constrainedl}, so a global minimizer
of~\eqref{rk:QUBO} that happens to be feasible is optimal for the
constrained problem.  At the boundary the encoding thus keeps the
optimal value and all constrained minimizers; what can be lost is only
the uniqueness of the cardinality, and only at the neighboring
cardinalities $k\pm 1$ for which the respective bound is attained.

\section{Quadratic Penalties and Inner Solvers}
\label{sec:penalties_solvers}

In this Section we address the negative effects of choosing too high
quadratic penalties to inner solvers, providing arguments  supporting
the broadly accepted intuitions for the case of certain inner solvers.

\subsection{Landscape Roughening by the Quadratic Penalty}
\label{sec:landscape_roughening}

We start with an observation that the quadratic penalty erects
energy barriers between feasible solutions close to each other in
Hamming-distance, that, for high enough $\eta$-s grow linearly
with~$\eta$, making the QUBO landscape progressively rougher.  Let us
restrict ourselves to quadratic penalties here, setting $\lambda = 0$:
\begin{equation}\label{eq:pure_quad}
  \min_{x \in \{0,1\}^n}\;
    \tfrac{1}{2} x^\top Q\, x \;+\; \eta\,(\mathbf{e}^\top x - k)^2,
\end{equation}
which is the most common approach in the QUBO literature
(cf.~\cite{Glover_2022,Lucas_2014}).

\begin{obs}[Penalty-induced barriers]\label{obs:barriers}
Let $x$ and $x'$ be two feasible solutions ($\mathbf{e}^\top x =
\mathbf{e}^\top x' = k$) that differ in exactly one pair of
coordinates: $x'_i = 1 - x_i$ and $x'_j = 1 - x_j$ for some
$i \neq j$, with $x'_l = x_l$ otherwise (a swap move).  Any
single-flip path from $x$ to $x'$ must pass through an intermediate
state $\bar{x}$ with $\mathbf{e}^\top \bar{x} = k \pm 1$.  The
penalty contribution at $\bar{x}$ is
\[
  \eta\,(\mathbf{e}^\top \bar{x} - k)^2 \;=\; \eta,
\]
which is absent at both $x$ and $x'$.  The energy barrier for
traversing this intermediate state therefore includes an additive
contribution $\eta$. This contribution is present regardless of the
objective values at $x$, $\bar{x}$, and $x'$.
\end{obs}

This observation has direct consequences to local-search heuristic
exploring the solution space via single-variable flips, such as the
Metropolis-Hastings steps of simulated annealing. Such solvers'
dynamics can only move between feasible solutions by passing through
infeasible intermediate states as described in
Observation~\ref{obs:barriers}.  In simulated annealing the
acceptance probability of an uphill move $\Delta E > 0$ is
$\exp(-\Delta E / T)$.  At the intermediate state $\bar{x}$, the total
energy change includes both the penalty contribution~$\eta$ and the
objective change $\Delta q$ (an instance-dependent constant).  The
acceptance probability is therefore bounded above by
$\exp(-(\eta - |\Delta q|) / T)$: once $\eta$ exceeds the
objective-change scale, the barrier is dominated by the penalty and
the escape probability vanishes exponentially in~$\eta$ at any fixed
temperature.

Meanwhile, the presented argument provides a rough intuition for
quantum annealing, too: --- larger $\eta$ implies higher barriers and
harder tunneling or thermal hopping between feasible solutions. The
actual mechanism is more complicated though: a rigorous analysis would
require tracking the minimum spectral gap over the full anneal
schedule, which is beyond the scope of this work. In this case,
however, another mechanism degrading solver performance for big
$\eta$-s can be identified, which we discuss next.

\subsection{Quantum Annealing: Hardware Precision}
\label{sec:distortion_qa}

At the beginning of the discussion we assume the problem is trivially
embedded, so that the realized physical qubits are the same as the
logical qubits belonging to the problem variables. (We will return to
the question of embedding later.) Thus a QUBO of the form
$f_{\lambda,\eta}(x) = \tfrac{1}{2}x^\top M(\lambda,\eta)\,x$ is
converted to an Ising Hamiltonian
\[
  H(s) \;=\; -\!\sum_{i<j} J_{ij}\,s_i s_j
             \;-\; \sum_i h_i\,s_i
             \;+\; \mathrm{const},
  \qquad s_i \in \{-1,+1\},
\]
via the standard substitution $x_i = (1+s_i)/2$, which maps each
off-diagonal QUBO entry to an Ising coupling of strength
$J_{ij} = M_{ij}(\lambda,\eta)/4$ and each diagonal entry to a
combination of local field~$h_i$ and constant offset.  The Ising
Hamiltonian is then programmed onto the physical device. The
analog control circuitry of the quantum processing unit (QPU) accepts coupler strengths $J_{ij}$ and local
fields $h_i$ only within a bounded dynamic range, and the programmed
values are subject to a finite precision floor~$\delta$: differences
smaller than $\delta$ are indistinguishable to the hardware.  We will
show that as the consequence of these two limitations, a large
quadratic penalty~$\eta$ gradually suppresses the objective's
contribution to the programmed Hamiltonian.

Let us investigate how $\eta$ rescales the dynamic range.
Let $J_{\max}^{\text{device}}$ denote the maximum absolute coupling
strength supported by the device.  To be acceptable for the hardware,
a QUBO whose Ising equivalent has maximum absolute coupling
$\|J\|_{\max}^{\text{problem}} := \max_{i,j}|J_{ij}|$ must be
uniformly rescaled by a factor
\[
  s \;=\; \frac{\|J\|_{\max}^{\text{problem}}}{J_{\max}^{\text{device}}},
\]
so that every programmed coupler's strength becomes $J_{ij}/s$.  The rescaling is
global: it applies identically to every entry, whether it encodes
objective information or penalty structure.

Now we show how the penalty term gradually dominates the objective
when increasing $\eta$.  Consider the pure quadratic penalty
encoding~\eqref{eq:pure_quad} ($\lambda = 0$), i.e., no linear penalty term.  The off-diagonal
entries of the QUBO matrix $M(0,\eta)$ take the form $\eta + c_{ij}$,
where $c_{ij}$ carries the objective information.  The QUBO-to-Ising
conversion scales all off-diagonal entries by the same factor
of~$1/4$, so the ratio between any two coupling strengths is preserved.  The
largest Ising coupling before rescaling is
$\|J\|_{\max}^{\text{problem}} = (\eta + \max c_{ij})/4$.
Let $J_{ij}, J_{i'j'}$
be the coupling strengths realizing $\max|c_{ij} - c_{i'j'}|$; after rescaling
by~$s$, their difference becomes
\begin{equation}\label{eq:coupler_diff}
  \Delta J_{\mathrm{scaled}}
  \;:=\; \left| \frac{J_{ij}}{s} - \frac{J_{i'j'}}{s} \right|
  \;=\; \frac{\max|c_{ij} - c_{i'j'}| \cdot J_{\max}^{\text{device}}}
             {\eta + \max(c_{ij})},
\end{equation}
where the factor of~$1/4$ from the Ising conversion cancels between
numerator and denominator.  When this quantity falls below the
hardware precision floor~$\delta$, the Hamiltonian the QPU actually
anneals is (to within its own precision) indistinguishable from a
graph with no edge-weight structure at all, and only the
cardinality penalty remains ``visible'' to the hardware.
The threshold is
\begin{equation}\label{eq:eta_threshold}
  \eta_{\text{device}}
  \;=\; \frac{\max|c_{ij} - c_{i'j'}| \cdot J_{\max}^{\text{device}}}{\delta}
    \;-\; \max(c_{ij}).
\end{equation}

An analogous compression affects the local fields~$h_i$: the
quadratic penalty contributes a uniform offset to every bias,
washing out the vertex-level contribution of the original objective after rescaling.
The wider bias range of QPU hardware (e.g.,
$h_{\max} = 6$ vs.\ $J_{\max}^{\text{device}} = 1$ on Advantage2) provides
partial protection, but the offset grows with both~$\eta$ and the
problem size, so it eventually dominates.  As shown in
Section~\ref{sec:sks_hw_threshold} for the SkS problem, a nonzero
linear penalty~$\lambda$ can partially or fully cancel this
offset, providing a second, bias-level argument for
preferring~$\lambda$ over~$\eta$.

Now let us address embedding.  On devices with limited connectivity
each logical variable is represented by a chain of physical qubits
coupled ferromagnetically with chain strength $J_{\mathrm{chain}}$.
The chain coupling is itself a penalty: it must be large enough that
breaking a chain never pays off, otherwise the embedded Ising problem
is not equivalent to the logical one.  A standard sufficient condition
for variable~$i$ is
\[
  |J_{\mathrm{chain}}| \;\ge\; |h_i| + \sum_{j}|J_{ij}|
\]
\cite{Choi_2008}.  For a penalised QUBO with a cardinality constraint
every pair of variables is coupled, so the right-hand side grows like
$(n-1)\eta/4$ plus the $\lambda$- and $\eta$-dependent local field; it
is much larger than $\|J\|_{\max}^{\text{problem}}$.  The chain
couplings share the hardware range with the problem couplings, so the
rescaling factor becomes
\[
  s \;=\; \max\Bigl(\frac{\|J\|_{\max}^{\text{problem}}}{J_{\max}^{\text{device}}},\;
          \frac{\|h\|_{\max}}{h_{\max}^{\text{device}}},\;
          \frac{|J_{\mathrm{chain}}|}{J_{\max}^{\text{device}}}\Bigr).
\]
As a consequence, the choice of chain strength introduces complications
similar to that of~$\eta$. Namely, a value below the bound lets chains break,
and the device samples a different problem. Values at the bound
set the rescaling factor $s$ globally, and thus compress every objective coupling 
by a factor growing with $n\eta$. 
The quadratic penalty thus degrades the embedded problem
twice: directly through~$\|J\|_{\max}^{\text{problem}}$ and indirectly
through the chain strength it requires.  Higher-connectivity
topologies (such as Zephyr on Advantage2, with maximum degree~20 vs.\
Pegasus's~15) yield shorter chains, mitigating but not eliminating the
effect.  Djidjev~\cite{Djidjev_2023b} treats the chain couplings
themselves by an augmented Lagrangian method; a joint treatment of the
problem and chain penalties is left for future work.

Let us now discuss the effect of the linear penalty on hardware
precision.  The linear penalty term $-\lambda(\mathbf{e}^\top x - k)$
contributes only to the diagonal entries of the QUBO matrix,
which map to local fields~$h_i$ in the Ising formulation.  The
off-diagonal elements~$J_{ij}$ --- where the pairwise objective
information resides --- are entirely unaffected by~$\lambda$:
increasing~$\lambda$ while decreasing~$\eta$ strictly increases
$\Delta J_{\mathrm{scaled}}$.

The effect on the local fields is problem-dependent.  The quadratic
penalty contributes a uniform bias offset that grows
with~$\eta$; the linear penalty contributes an offset of opposite
sign.  Depending on the problem, these two offsets can partially
or fully cancel (see Section~\ref{sec:sks_hw_threshold} for the
SkS case), so decreasing~$\eta$ while increasing~$\lambda$ can
simultaneously increase $\Delta J_{\mathrm{scaled}}$ and reduce
the penalty-induced offset in the local fields.
The wider bias range of the Advantage2 system
($\texttt{h\_range} = [-6,\,6]$~\cite{Advantage2WP}, up from
$[-4,\,4]$ of the Advantage) provides additional headroom for
accommodating whichever offset remains after cancellation.

In summary, on a fully connected device $\lambda$ is neutral on
couplings and potentially beneficial on biases, while $\eta$ is
harmful on both.  On embedded hardware, $\lambda$ is not free: it
enlarges $|h_i|$ and therefore the chain strength required for an
equivalent embedding, although far less than $\eta$, which enters
every coupler of the chain bound.  This provides a hardware-grounded
motivation, in addition to the landscape arguments of
Section~\ref{sec:landscape_roughening}, for keeping both penalties
small and preferring~$\lambda$ over~$\eta$ whenever the problem
structure permits (cf.\ Corollary~\ref{rk:lambda}).  In
Section~\ref{sec:sks_hw_threshold} we evaluate the coupling threshold
for the SkS problem.

The precision argument above concerns the static programming of
the Ising Hamiltonian and is a necessary condition only: once the
objective's contribution falls below the precision floor, no anneal
schedule can recover it.  It does not by itself predict how a QPU
performs above that floor, which also depends on chain integrity,
noise and the minimum spectral gap along the anneal; a rigorous
treatment of the latter is beyond the scope of this work.  In our
experiments (Section~\ref{sec:cs_results}) these factors cannot be
separated from the landscape effect of
Section~\ref{sec:landscape_roughening}: 
the effects are also present when steepest descent from
random starts is used (Section~\ref{sec:cs_solvers}), without any
quantum hardware.

\section{Application to the Sparsest \texorpdfstring{$k$}{k}-Subgraph Problem}
\label{sec:application_sks}

In order to conduct a more detailed illustrative analysis of the
importance of the right choice of penalty coefficients, we now
specialize further by setting $Q := A$, the adjacency matrix of a
simple undirected graph $G$ on $n$ vertices.  As the graph is not
weighted and undirected,
\begin{equation*}
  A_{ij} = A_{ji}\in \{0,1\},\  A_{ii} = 0 \quad \quad \forall i,j,
\end{equation*}
thus, compared to the problem in
Section~\ref{sec:card_qbp}, we have a more specific form of the
objective function. (The symmetry of $Q$ has been assumed
wlog. also thus far.)  With these assumptions, our problem becomes a
constrained quadratic binary formulation of the well-known NP-hard
Sparsest $k$-Subgraph problem (SkS): find an induced subgraph
of cardinality $k$ of $G$ with as few edges as possible.  Any
minimizer $x_k^*$ of~\eqref{rk:QBO_constrainedl} with $Q := A$ induces
a $k$-subgraph $G[x_k^*]$ with $m_k = \bqp{A}{x_k^*}$ edges --- the
minimum achievable.

\subsection{Complexity of Finding the Optimal $\eta$}

Theorem~\ref{thm:min_eta} characterizes the minimum penalty
$\eta^*$ in terms of $A_k, B_k$ and their unit-step counterparts.
We now show that computing these quantities is intractable in
general, even in the SkS setting ($Q := A$).

\begin{thm}\label{thm:np_hard}
Computing the values $(A_k^1-B_k^1)/2$ and $(A_k-B_k)/2$ is
$NP$-hard.
\end{thm}

\begin{proof}
We reduce from the maximum independent set problem.  Let $G$ be an
undirected graph on~$n$ vertices with adjacency matrix~$A$, and let
$\alpha(G)$ denote the independence number of~$G$.

For any $k \leq \alpha(G)$, every $k$-element subset can be chosen
to be independent, so $q_l = 0$ for all $l \leq \alpha(G)$;
consequently $A_k^1 = B_k^1 = A_k = B_k = 0$ whenever
$k < \alpha(G)$.

At $k = \alpha(G)$: we still have $q_k = q_{k-1} = 0$, so
$A_k^1 = A_k = 0$.  However, $q_{k+1} > 0$ since no independent
set of size $k+1$ exists, giving $B_k^1 = q_{k+1} > 0$ and hence
$(A_k^1 - B_k^1)/2 < 0$.

Therefore the first~$k$ for which $(A_k^1 - B_k^1)/2$ is nonzero
equals $\alpha(G)$.  An oracle that computes $(A_k^1 - B_k^1)/2$
for any given~$k$ would determine $\alpha(G)$ via binary search,
and since computing $\alpha(G)$ is $NP$-hard, so is computing
$(A_k^1 - B_k^1)/2$.  The argument for $(A_k - B_k)/2$ is
identical. \qed
\end{proof}
This has an immediate consequence relevant for the practitioner:
\begin{cor}\label{cor:np_hard_eta}
Computing the optimal (smallest) value of the quadratic penalty
coefficient for the QUBO encoding~\eqref{rk:QUBO} of the SkS
problem is $NP$-hard.
\end{cor}
\begin{proof}
By Theorem~\ref{thm:min_eta}, the infimum $\eta^*$ is determined
by $A_k, B_k, A_k^1, B_k^1$.  Since computing these quantities is
$NP$-hard (Theorem~\ref{thm:np_hard}), so is computing~$\eta^*$. \qed
\end{proof}

Theorem~\ref{thm:np_hard} shows that finding $\eta^*$ exactly is
as hard as the original problem.  Nevertheless, the interval
$[\tfrac{A_k^1-B_k^1}{2},\,\tfrac{A_k-B_k}{2}]$ of
Theorem~\ref{thm:min_eta} can collapse to a single point, making
the bounds tight:
\begin{cor}\label{cor:tight_lb}
If\/ $A_k = A_k^1$ and $B_k = B_k^1$, then
$\eta^* = (A_k^1-B_k^1)/2 = (A_k-B_k)/2$: the lower and upper
bounds on the infimum coincide, and equivalence holds for any
$\eta > \eta^*$.
\end{cor}
\begin{proof}
Substituting $A_k = A_k^1$ and $B_k = B_k^1$ into
Theorem~\ref{thm:min_eta}(b) gives
$\eta^* \in [(A_k^1-B_k^1)/2,\,(A_k^1-B_k^1)/2]$, a single
point. \qed
\end{proof}

The tightness condition is non-vacuous.  For the path graph~$P_n$,
the optimal value function $q_l$ (minimum number of edges in an
induced $l$-subgraph) is convex in~$l$: each additional vertex
contributes at most one new edge, and this marginal cost is
non-decreasing.  Convexity of $q_l$ implies that the difference
quotients $(q_k - q_l)/(k-l)$ are non-decreasing in the step size
$k-l$ and that $(q_l - q_k)/(l-k)$ are non-increasing, so the
extrema in the definitions of $A_k$ and $B_k$ are achieved at
$l = k \mp 1$.  Hence $A_k = A_k^1$ and $B_k = B_k^1$ for
all~$k$.

\subsection{Calculable Penalty Bounds}

Given the $NP$-hardness of finding the optimal $\eta$ in general,
to encode an actual problem instance one can use estimates of $A_k$
and $B_k$ from~\eqref{rk:ABvals} to compute practical penalty
coefficients based on Theorem~\ref{rk:theorem}.  The following
Corollary provides computable values of $\lambda$ and $\eta$ for
the sparsest $k$-subgraph problem.  Note that $\valapprox{m}$ can
be efficiently calculated using any heuristic, such as the greedy
algorithm that starts from the full graph~$G$ and iteratively
removes the vertex of maximum degree until a $k$-subgraph remains
(as described in Section~\ref{augmented_lagrangian} and formalized
in Observation~\ref{obs:greedy_init}).

\begin{cor}[Sufficiency bounds]\label{cor:penalty_bounds}
Let $Q := A$ be the adjacency matrix of a graph $G$ with $n$ vertices, and let
$\valapprox{m} := \bqp{A}{x}$ for some $\inbin{x}$ with $\cardy{x} = k$.
Let $x^*$ be a minimizer of the QUBO with parameters $\lambda$ and $\eta$.

\begin{equation}
    \eta > \min\!\left\{\dfrac{\valapprox{m}}{2},\,\dfrac{k-1}{2}\right\}
    \quad and \quad
    \lambda = \min\!\left\{\dfrac{\valapprox{m}}{2},\,\dfrac{k-1}{2}\right\}\tag{T1}      
\end{equation}
\begin{equation}
    \eta > \min\!\left\{\valapprox{m},\,k-1\right\}
    \quad and \quad
    \lambda = 0\tag{T2}
\end{equation}

Then $G[x^*]$ is the sparsest $k$-subgraph of $G$.
\end{cor}
\begin{proof}
  Both items follow from Theorem~\ref{rk:theorem} together with the
  bounds $\valapprox{m} \geq m_k \geq A_k$, $k-1 \geq A_k$ and $0\leq B_k$.\\[4pt]
  \emph{(T1):} Set
  $\valapprox{A} := \min\{\valapprox{m},\,k-1\}$ and
  $\valapprox{B} := 0$.  Since $\valapprox{A} \geq A_k$ and $\valapprox{B} \leq B_k$, the condition
  of Theorem~\ref{rk:theorem} is satisfied with
  $\lambda = \min\{\valapprox{m}/2,\,(k-1)/2\}$ and $\eta = 
  \valapprox{A}/2+\epsilon, \ \epsilon > 0$.\\[4pt]
  \emph{(T2):} With the same $\valapprox{A}$ and $\valapprox{B}$,
  choose $\epsilon' := \valapprox{A}/2$, so that
  $\lambda = \valapprox{A}/2 - \epsilon' = 0$ is admissible, and any
  $\epsilon > \epsilon'$.  The condition of Theorem~\ref{rk:theorem} is
  then satisfied with $\eta = \valapprox{A}/2 + \epsilon$, that is, with
  any $\eta > \valapprox{A}$.
  
  \qed
\end{proof}
\noindent
Both bounds are computable from any feasible heuristic solution  without solving any auxiliary optimization problem.

The inequalities for $\eta$ in Corollary~\ref{cor:penalty_bounds} are
strict.  If $\min\{\valapprox{m},\,k-1\} = 0$, then (T1) and (T2)
coincide: $\lambda = 0$ and any $\eta > 0$ is admissible.  Otherwise
the boundary values themselves, which are the infima of the admissible
quadratic penalties, are covered by Proposition~\ref{prop:boundary}.
Let $m_l := q_l$ denote the number of edges of a sparsest $l$-subgraph
of~$G$.

\begin{cor}[Boundary values]\label{cor:sks_boundary}
Let $Q := A$ and $m_k > 0$, and let $\lambda$ and $\eta$ be given by
(T1) or by (T2) of Corollary~\ref{cor:penalty_bounds} with the strict
inequality for $\eta$ replaced by equality.  Then the optimal value of
the QUBO~\eqref{rk:QUBO} is $m_k$ and every sparsest $k$-subgraph is a
minimizer.  All minimizers have cardinality~$k$, unless
\begin{equation}\label{eq:sks_tie}
  m_k - m_{k-1} \;=\; \min\{\valapprox{m},\,k-1\},
\end{equation}
in which case the minimizers are exactly the sparsest $k$-subgraphs
and the sparsest $(k-1)$-subgraphs.  Condition~\eqref{eq:sks_tie} can
hold only if $m_{k-1} = 0$ (that is, if $k-1$ is the independence
number of~$G$) or if the complement of~$G$ is disconnected.
\end{cor}
\begin{proof}
As $m_1 = 0$, we have $k \geq 2$, and $\valapprox{m} \geq m_k > 0$, so
$\eta > 0$.  In both cases
$\lambda+\eta = \min\{\valapprox{m},\,k-1\} \geq \min\{m_k,\,k-1\}
\geq A_k$ and $\lambda-\eta \leq 0$.  For $l>k$ we have
$m_l \geq m_k > 0$, so a sparsest $l$-subgraph has a vertex of
positive degree, and removing it shows that $m_{l-1} \leq m_l - 1$.  Hence
$m_l - m_k \geq l-k$ for all $l>k$, so that $B_k \geq 1$ and
$B_k^1 \geq 1 > \lambda-\eta$.  Proposition~\ref{prop:boundary} now
gives the optimal value, excludes minimizers of cardinality $k+1$, and
shows that minimizers of cardinality $k-1$ exist if and only if
$A_k^1 = \lambda+\eta$, which is~\eqref{eq:sks_tie}; these are the
vectors of cardinality $k-1$ with objective value
$m_k - (\lambda+\eta) = m_{k-1}$.

If~\eqref{eq:sks_tie} holds, then
$m_k - m_{k-1} \leq \min\{m_k,\,k-1\} \leq \min\{\valapprox{m},\,k-1\}$
holds with equality throughout, so $m_k - m_{k-1} \in \{m_k,\,k-1\}$.
In the first case $m_{k-1} = 0$.  In the second case let $S$ be the
vertex set of a sparsest $(k-1)$-subgraph.  For every vertex
$v \notin S$ the subgraph induced by $S \cup \{v\}$ has
$m_{k-1} + \deg_S(v) \geq m_k = m_{k-1} + k-1$ edges, where
$\deg_S(v)$ is the number of neighbors of $v$ in~$S$; hence $v$ is
adjacent to all vertices of~$S$.  The complement of~$G$ then has no
edge between $S$ and its complement, both of which are non-empty. \qed
\end{proof}

\subsection{Quantum Hardware Precision Threshold for the SkS Problem}
\label{sec:sks_hw_threshold}

We specialize the coupling-difference and threshold formulas of
Section~\ref{sec:distortion_qa} for the SkS problem at $k = n/2$,
where the penalty offset in the local fields vanishes
($|n - 2k| = 0$) and the coupling threshold is the only limitation.

The off-diagonal QUBO entries read
$M_{ij} = \tfrac{1}{2}A_{ij} + \eta$, with the maximal difference in
the objective coefficients $\max|c_{ij} - c_{i'j'}| = 1/2$ and
$\max c_{ij} = 1/2$.
Substituting into~\eqref{eq:coupler_diff}:
\begin{equation}\label{eq:sks_coupler_signal}
  \Delta J_{\mathrm{scaled}}
  \;=\; \frac{J_{\max}^{\text{device}}}{2\eta + 1}.
\end{equation}
The graph structure becomes invisible to the hardware when
$\Delta J_{\mathrm{scaled}} < \delta$; substituting into~\eqref{eq:eta_threshold}
gives the threshold:
\begin{equation}\label{eq:eta_threshold_sks}
  \eta \;>\; \frac{J_{\max}^{\text{device}}}{2\delta} - \frac{1}{2}
  \;=:\; \eta_{\text{device}}.
\end{equation}

These are the values for a trivially embedded problem.  On an embedded
device the rescaling factor is set by the larger of the problem and
chain scales (Section~\ref{sec:distortion_qa}), so
$\Delta J_{\mathrm{scaled}}$ can be much smaller than
\eqref{eq:sks_coupler_signal}; with a chain strength that exceeds every
problem coupling it no longer depends on~$\eta$ at all
(Section~\ref{sec:cs_results}).

\begin{bem}
For $k \neq n/2$ the quadratic penalty also compresses the local
fields~$h_i$, introducing a second, potentially tighter threshold
that depends on~$|n - 2k|$ and can be partially canceled by a
nonzero~$\lambda$.  A full analysis of the interplay between the
coupling and bias thresholds across the entire $(n, k, \rho)$
parameter space is the subject of ongoing work.
\end{bem}

\section{Computational Study: Setup}
\label{sec:computational_study}

Having discussed the theoretical aspects of our contribution, let us
now turn our attention to numerical experimentation. This Section is
restricted to SkS problems, where the QUBO matrix $Q$ equals to an
adjacency matrix $A$ of a simple undirected graph, thus it is a
symmetric matrix with zeros in its diagonal, and $0$-s or $1$-s in the
off-diagonal elements.

\subsection{Practical Augmented Lagrangian Method}
\label{augmented_lagrangian}

We now introduce a practical algorithm that adjusts the penalty
parameters $(\lambda, \eta)$ of the QUBO
encoding~\eqref{rk:QUBO} of a cardinality-constrained
problem~\eqref{rk:QBO_constrainedl} until the cardinality
constraint is satisfied.  The starting point is the observation
that the optimal penalties of Theorem~\ref{thm:min_eta} can be
approximated from a single greedy pass.

Let us therefore describe the greedy algorithm.
Given a graph $G$ on $n$ vertices and a target cardinality~$k$, the
greedy algorithm for the sparsest $k$-subgraph proceeds by
successive vertex removal. Starting from $G_n := G$ as the full graph, it
iteratively removes the vertex of maximum degree in the current
induced subgraph, producing a nested sequence
$G_n \supset G_{n-1} \supset \cdots \supset G_k$.  At
each step, the number of edges decreases by the degree of the
removed vertex:
\begin{equation}\label{eq:greedy_edge_drop}
  g_l - g_{l-1} \;=\; \Delta(G_l),
\end{equation}
where $g_l$ denotes the edge count of the greedy $l$-subgraph and
$\Delta(G_l)$ the maximum degree in $G_l$.  The subgraph
$G_k$ serves both as a feasible heuristic solution (providing
the edge-count estimate $\tilde{m}_k := g_k$ used in
Corollary~\ref{cor:penalty_bounds}) and as the source of the
penalty initialization derived below.

\begin{obs}[Greedy initialization of $\lambda$ and $\eta$]
\label{obs:greedy_init}
For the penalties $\lambda^* = (A_k^1+B_k^1)/2$ and
$\eta^* = (A_k^1-B_k^1)/2$ of Theorem~\ref{thm:min_eta}:
\begin{enumerate}
  \item[(i)] The unit-step quantities
    (Definition~\ref{def:unit_step}) satisfy
    \[
      A_k^1 \;\approx\; \Delta(G[g_k]),
      \qquad
      B_k^1 \;\approx\; \Delta(G[g_{k+1}]).
    \]
  \item[(ii)] If the maximum degree in the greedy subgraph varies
    slowly with cardinality, i.e.,
    $\Delta(G[g_k]) \approx \Delta(G[g_{k+1}])$, then
    \[
      \lambda^* \;\approx\; \Delta(G[g_k]),
      \qquad
      \eta^* \;\approx\; 0.
    \]
\end{enumerate}
\end{obs}

\begin{proof}
For part~(i), write $m_l^* = g_l - \epsilon_l$ where
$\epsilon_l \geq 0$ is the suboptimality gap of the greedy
$l$-subgraph.  By Definition~\ref{def:unit_step},
\[
  A_k^1 = m_k^* - m_{k-1}^*
    = (g_k - g_{k-1}) - (\epsilon_k - \epsilon_{k-1}).
\]
From~\eqref{eq:greedy_edge_drop}, $g_k - g_{k-1} =
\Delta(G[g_k])$: this is the degree of the vertex removed in the
greedy step from cardinality~$k$ to $k{-}1$.  When the
suboptimality gap varies slowly with cardinality
($\epsilon_k - \epsilon_{k-1}$ small), the claim
$A_k^1 \approx \Delta(G[g_k])$ follows.  The argument for
$B_k^1$ is analogous: $g_{k+1} - g_k = \Delta(G[g_{k+1}])$ is
the degree of the vertex whose removal produced $G[g_k]$ from
$G[g_{k+1}]$.

Part~(ii) follows immediately from part~(i) by substituting
$\Delta(G[g_k]) \approx \Delta(G[g_{k+1}])$ into the definitions
of $\lambda^*$ and $\eta^*$. \qed
\end{proof}

The approximation $\Delta(G[g_k])$ is computable as a by-product
of the greedy pass at negligible additional cost.  Since
$\epsilon_k - \epsilon_{k-1}$ can be of either sign,
$\Delta(G[g_k])$ is neither a guaranteed upper nor lower bound for
$\lambda^*$ --- it is a practical estimate.

Part~(ii) has an important consequence: for graphs where the
maximum degree in the greedy subgraph is approximately constant
near cardinality~$k$, the optimal quadratic penalty is close to
zero and the constraint can be enforced almost entirely through
the linear term~$\lambda$.  This motivates
Algorithm~\ref{alg:augmented_lagrangian}, which warm-starts
$\lambda_0 \approx \Delta(G[g_k])$ and $\eta_0$ at a small value,
then adjusts both iteratively.

\begin{algorithm}[htbp]
\caption{Augmented Lagrangian method for the SkS problem}
\label{alg:augmented_lagrangian}
\begin{algorithmic}[1]
\Require Initial $\lambda \ge 0$, $\eta > 0$,
  $k \in \{1,2,\dots,n-1\}$, penalty scaling factor $\beta > 1$
\Repeat
    \State $x \gets \arg\min_{x \in \{0,1\}^n}\;
      \tfrac{1}{2}x^\top\! A\,x
      - \lambda(\mathbf{e}^\top\! x - k)
      + \eta(\mathbf{e}^\top\! x - k)^2$
    \If{$\mathbf{e}^\top\! x \neq k$}
       \State $\lambda \gets
         \lambda + 2\eta\,(k - \mathbf{e}^\top\! x)$
         \Comment{standard augmented Lagrangian multiplier update}
    \EndIf
    \State $\eta \gets \beta\,\eta$
\Until{$\mathbf{e}^\top\! x = k$ or iteration limit reached}
\end{algorithmic}
\end{algorithm}

The algorithm implements a standard augmented Lagrangian method:
at each iteration, $\lambda$ is corrected by the signed
cardinality violation scaled by the current penalty $2\eta$,
and $\eta$ is increased by a factor $\beta > 1$.  The update
$\lambda \leftarrow \lambda + 2\eta(k - \mathbf{e}^\top x)$ is
the classical multiplier update of
Hestenes~\cite{Hestenes_1969} and Powell~\cite{Powell_1969},
specialized to a single equality constraint; it ties the step
size to the current penalty level, so that early iterations
(small~$\eta$) make conservative $\lambda$ adjustments while
later iterations (larger~$\eta$) correct more aggressively.

We note here that an oscillation of $\lambda$ can arise, which can be
circumvented by damping. Here we describe this phenomenon and the
damping strategy which will be termed as ``GLO variant'' in what
follows.  Algorithm~\ref{alg:augmented_lagrangian} is designed for
continuous problems where the constraint violation decreases
monotonically toward zero.  For the SkS problem the inner solver
returns a discrete cardinality, and it is common for the returned
value to alternate between $k{-}1$ and $k{+}1$ without ever landing
on~$k$.  Because the violation flips sign while $\eta$ grows, the
$\lambda$ corrections grow and alternate in sign, producing a
diverging oscillation.

The oscillation can be suppressed by observing that two consecutive
alternations bracket the correct multiplier~$\lambda^*$.  Let
$v_i = k - \mathbf{e}^\top x_i$ denote the signed violation at
iteration~$i$.  Once two consecutive sign changes are detected
(i.e., $v_{i-1}$ and $v_i$ have opposite signs for two successive
iteration pairs), $\lambda$ is locked to
\begin{equation}\label{eq:lambda_lock}
  \lambda^* \;=\;
    \frac{|v_{i-1}|\,\tilde{\lambda}_i
          \;+\; |v_i|\,\tilde{\lambda}_{i-1}}
         {|v_{i-1}| + |v_i|},
\end{equation}
where $\tilde{\lambda}_j = \lambda_j + 2\eta_j\, v_j$ is the shadow
multiplier (the value $\lambda$ \emph{would} take after the standard
update at iteration~$j$).  For the common symmetric case
$|v_{i-1}| = |v_i| = 1$ this reduces to the simple average
$\lambda^* = (\tilde{\lambda}_i + \tilde{\lambda}_{i-1})/2$.  After
locking, $\lambda$ is no longer updated; only $\eta$ continues to
grow.  In our experiments, the locked $\lambda^*$ typically results in
$\mathbf{e}^\top x = k$ on the next use of the inner solver.

\subsection{Penalty Strategies}
\label{sec:cs_variants}

The experiments compare the following three fundamentally different
strategies for choosing the penalty coefficients $(\lambda,\eta)$ in
Algorithm~\ref{alg:augmented_lagrangian}, each motivated by a
different level of problem-specific information.

\paragraph{AL (Augmented Lagrangian with Greedy warm-start).}
The linear multiplier $\lambda_0$ is warm-started from the maximum
degree of the greedy $k$-subgraph (Observation~\ref{obs:greedy_init}),
while the quadratic penalty starts at a small fixed value $\eta_0$.
The AL loop then adjusts both $\lambda$ and $\eta$ iteratively until
the cardinality constraint is satisfied.  The GLO variant adds the
oscillation damping described in Section~\ref{augmented_lagrangian},
preventing the divergent $\lambda$ oscillation that can occur with
discrete cardinality returns.

\paragraph{T1 (Corollary~\ref{cor:penalty_bounds}, symmetric split).}
Both penalties are set from the Corollary~\ref{cor:penalty_bounds}:
$\lambda_0 =, \eta_0 \ge \min\{\tilde{m}_k/2,\,(k{-}1)/2\}$, where
$\tilde{m}_k$ is an estimate of the optimal edge count~$m_k^*$. The values are
instance-dependent; 
Corollary~\ref{cor:penalty_bounds} guarantees that every global QUBO
minimizer has cardinality~$k$ provided that $\tilde{m}_k \geq m_k^*$ and $\eta_0 > \min\{\tilde{m}_k/2,\,(k{-}1)/2\}$. 

\paragraph{T2 (Corollary~\ref{cor:penalty_bounds}, pure quadratic).}
We set $\eta_0\ge \min\{\tilde{m}_k,\,k{-}1\}$, $\lambda_0 = 0$,
implementing the commonly used quadratic penalty encoding. This can
serve as a reference point for measuring the benefit of a nonzero
linear term $\lambda$.  Note that by construction
$\eta_{T2} = 2\eta_{T1}$ instance-by-instance. 
The experiments of Section~\ref{sec:computational_results} use the
known optimum, $\tilde{m}_k = m_k^*$; the comparison with other
penalty bounds in Section~\ref{sec:cs_cor} uses the greedy value
$\tilde{m}_k = g_k$.

\paragraph{Boundary values.}
In the experiments of
Section~\ref{sec:computational_results}, T1 and T2 are run with the
values given above, that
is, with equality in place of the strict inequality for $\eta$ in
Corollary~\ref{cor:penalty_bounds}.  These are the infima of the
quadratic penalties the Corollary~\ref{cor:penalty_bounds} admits: every admissible $\eta$ is
larger.  The only exception are the instances with
$\min\{\tilde{m}_k,\,k{-}1\} = 0$; with $\tilde{m}_k = m_k^*$ these are
the five instances with $n = 40$, $k = 10$ and $\bar{\rho} = 25\%$.
There $\lambda_0 = 0$ and $\eta_0$
is set to a small positive value, as
the Corollary requires.  By Corollary~\ref{cor:sks_boundary}, at the
boundary values the QUBO has optimal value $m_k^*$ and every sparsest
$k$-subgraph is among its global minimizers.  Compared with a strictly
larger~$\eta$, the only difference is that the sparsest
$(k{-}1)$-subgraphs are additional global minimizers on instances
satisfying~\eqref{eq:sks_tie}.

With $\tilde{m}_k = m_k^*$, condition~\eqref{eq:sks_tie} holds for none
of the 270 benchmark instances (Section~\ref{sec:cs_instances}).  All
original (uncomplemented) benchmark graphs are connected, which
excludes $m_k^* - m_{k-1}^* = k-1$; and on the 21 instances with
$1 \leq m_k^* \leq k-1$ an exhaustive search shows that $G$ has no
independent set of size $k-1$, so that $m_{k-1}^* > 0$.  Hence, on all
instances, the global minimizers of the T1 and T2 QUBOs used in
Section~\ref{sec:computational_results} are exactly the sparsest
$k$-subgraphs, as they would be for any strictly larger~$\eta$.

Running each of these strategies with several inner solvers serves one
of the central goals of our experimental agenda: to study how the
penalty strategy and the inner solver interact. The penalty strategy
determines the structure of the QUBO to be solved by the inner solver,
and different inner solvers can have different strengths and
weaknesses in solving it. Hence, it is instructive to understand the
extent to which each can compensate for shortcomings of the
other.

In addition to the three penalty strategies, we report the greedy
algorithm of Section~\ref{augmented_lagrangian} run on its own
(Greedy) as a penalty-free baseline.  It constructs a $k$-subgraph
directly by successive vertex removal, without a QUBO formulation or
an inner solver, so its output always satisfies the cardinality
constraint.  Its objective value is $g_k$, the same quantity that
warm-starts AL and serves as the estimate $\tilde{m}_k$ where $m_k^*$
is unknown (Observation~\ref{obs:greedy_init}).

\subsection{Inner Solvers}
\label{sec:cs_solvers}

Algorithm~\ref{alg:augmented_lagrangian} is solver-agnostic: any QUBO
solver can serve as the inner solver at line~2.  To test the
generality of the penalty-distortion effects across solver families,
we evaluate a broad range of solvers. These include steepest descent
from random starts as a baseline and four heuristic solvers.  Thus the
solver approaches cover classical local search, advanced
hyper-heuristics, physics-based simulation, and quantum hardware.

\paragraph{MQLib (classical hyper-heuristic).}
MQLib~\cite{MQLibRepo} is a C++ library of classical heuristics for
Max-Cut and QUBO, developed from the systematic study of Dunning,
Gupta, and Silberholz~\cite{MQLIB}. The version of the solver we are
using contains 39 different algorithms, covering numerous versions of
different types of heuristics (greedy, simulated
annealing, tabu search, etc.), collected from a systematic
bibliography search. We use its hyper-heuristic
mode, which selects a heuristic from the portfolio using a decision
tree based on calculable graph-theoretic features of the graph of the
MAX-CUT equivalent of the given QUBO instance.

MQLib was run on a consumer-grade computer equipped with an Intel
Core i7-1370P CPU at 2.2\,GHz, and 32\,GB RAM. This machine was used
for all the experiments.

\paragraph{Quantum Annealing (QPU).}
We use a D-Wave quantum annealer~\cite{Johnson_2011}, the Advantage2
system \texttt{Advantage2\_system2.1} (Zephyr topology, 4400+ qubits,
maximum degree~20, $\texttt{h\_range}=[-6,6]$,
$\texttt{j\_range}=[-1,1]$), accessed via the Leap cloud service and
the Ocean SDK~\cite{OceanSDK}.  Because every QUBO with $\eta>0$ is
dense (the quadratic penalty couples all pairs of variables), we use
one fixed minor embedding of the complete graph $K_n$ per problem
size, computed with \texttt{minorminer}, which implements the
heuristic of Cai, Macready, and Roy~\cite{Cai_2014}; mean chain
lengths grow from 4.6 ($n=40$) to 20.3 ($n=140$) physical qubits
(Table~\ref{tab:qpu_chains}).  Each annealing run draws 500 samples
(\texttt{num\_reads}) at an anneal time of $20\,\mu\text{s}$.  The
chain strength is fixed at~15 (Ising units) for all instances and all
encodings; this value was chosen by trial on a few instances, after
the default setting gave poor results.  The embedded Hamiltonian is
auto-scaled into the hardware ranges, and broken chains are resolved
by majority vote.  We evaluate two modes: raw QPU samples
(\textbf{QPU}), and QPU samples followed by a classical
steepest-descent post-processing step (\textbf{QPU+PP}), which
locally re-minimizes each sample using \texttt{SteepestDescentSolver}
from \texttt{dwave-samplers}~\cite{DwaveSamplers}.

\paragraph{Simulated Annealing (SA).}
We use D-Wave's
\texttt{SimulatedAnnealingSampler}~\cite{DwaveSamplers}, a
Metropolis-based local-search heuristic~\cite{Kirkpatrick_1983}
sharing the same data format as the QPU samplers. SA serves as the
standard classical counterpart of the QPU: it explores the same QUBO
landscape via thermal moves rather than quantum dynamics.  This solver
was run on the computer also used for the MQLib calculations.

\paragraph{Steepest descent from random starts (SD).}
As a classical baseline for the QPU+PP mode, we run the same
steepest-descent routine (\texttt{SteepestDescentSolver}) from
\texttt{num\_reads} uniformly random initial states instead of QPU
samples, inside the same penalty loop and on the same instances as
the QPU, with $\texttt{num\_reads}=500$ as for the QPU.  SD is a
deterministic single-flip descent with no noise, no embedding and no
precision limit.  Comparing QPU+PP with SD at equal
\texttt{num\_reads} isolates the contribution of the QPU samples as
starting points.  D-Wave's documentation uses the same post-processing
and the same random-start comparison to illustrate the benefit of
post-processing QPU samples~\cite{DWave_PP_Greedy}.

\paragraph{Simulated Bifurcation Machine (SBM).}
Simulated bifurcation has its roots in a quantum annealing proposal
based on coupled nonlinear quantum oscillators~\cite{Goto_2016}. It
solves the equations of motion of the classical (mechanical)
Hamiltonian system whose canonical quantization is the quantum
annealing proposal exhibiting bifurcations. It is a probabilistic QUBO
sampler, with randomness arising from the sensitivity to initial
conditions. The time step, the number of time steps and the number of
samples are hyperparameters whose settings can affect
performance. Beside the original continuous-time version, the
so-called ballistic and discrete-time
versions~\cite{Goto_2019,SBM,Kanao_2022} are also prevalent, and
recent modifications add further quantum-like
behavior~\cite{2604.01050}. SBM is able to handle large, dense QUBO
instances without the need of embedding.  We have used the cloud
implementation provided by Quantumz.io~\cite{QuantumzIO,VeloxQWiki}.

\subsection{Instances}
\label{sec:cs_instances}

In this study we use a benchmark set of 270 instances, derived
from the standard $k$-cluster benchmark of Ref.~\cite{bechmark}.
Because the $k$-cluster (densest $k$-subgraph) of a graph~$G$
corresponds to the sparsest $k$-subgraph of the complement~$\bar{G}$,
all instances are complemented before solving.  The benchmark
comprises graphs of six sizes
($n \in \{40, 80, 100, 120, 140, 160\}$), three edge densities
($\rho \in \{25\%, 50\%, 75\%\}$; the complemented graphs, on which
the SkS problem is solved, have density $\bar{\rho} = 1 - \rho$), three target cardinalities
($k \in \{n/4,\, n/2,\, 3n/4\}$), and five random replicates per
configuration, giving $6 \times 3 \times 3 \times 5 = 270$ instances
in total.

Not all 270 instances were solved by every algorithm variant.
QPU-based variants and their steepest-descent baseline (SD) were run
only on the 225 instances with $n \leq 140$; for the 45 instances with
$n = 160$, \texttt{minorminer} did not find an embedding of the
complete graph $K_{160}$ on the device.  SA, MQLib, SBM, and Greedy
were run on all 270 instances.

\section{Computational Study: Results}
\label{sec:computational_results}

We present results by solver family: quantum annealing with simulated
annealing as its classical counterpart (Section~\ref{sec:cs_results}),
MQLib and SBM (Section~\ref{sec:cs_mqlib}).
Where applicable, all three penalty strategies ---
AL, T1, and T2 --- are compared.

Throughout the experiments
(Sections~\ref{sec:cs_results}--\ref{sec:cs_mqlib}), the T1 and T2
penalties are calculated using the known optimum $m_k^*$ and taken at
their boundary values (Section~\ref{sec:cs_variants}), that is, at the
infimum of the quadratic penalties
Corollary~\ref{cor:penalty_bounds} admits.  If these variants perform
poorly, then even the tightest penalties the Corollary can provide are
too large for the solver in question.  Note that a positive error
means that the solver returned a feasible but suboptimal
$k$-subgraph, whose QUBO value is strictly above the minimum; such an
error cannot be caused by a tie with a $(k{-}1)$-subgraph
(Corollary~\ref{cor:sks_boundary}).  In
practice $m_k^*$ is unknown and has to be replaced by a heuristic
value such as $g_k$ (Observation~\ref{obs:greedy_init}), which can
only increase the T1 and T2 penalties.

\subsection{QPU, Steepest Descent and Simulated Annealing (D-Wave Ocean SDK)}
\label{sec:cs_results}

We first evaluate the QPU variants on the 225 instances with
$n\le 140$, together with two classical counterparts: steepest descent
from random starts (SD) on the same instances, and simulated annealing
(SA) on all 270 instances.  In line with the aim stated in
Section~\ref{sec:cs_variants}, this subsection varies both the penalty
strategy (AL, T1, T2) and the inner solver (QPU, SD, SA).  Under T1
and T2 the penalties depend only on the instance, so all solvers
receive the same QUBO; under AL all solvers start from the same
initial penalties, but each follows its own iteration path.  SD runs
the same post-processing routine as QPU+PP at the same
\texttt{num\_reads}, so the difference between these two isolates the
contribution of the QPU samples as starting points.

The QPU results are presented in Table~\ref{tab:sks_qpu_aggregate}.
The AL loop with post-processing (AL-QPU+PP) reaches 83.1\%
optimality with a mean error of 0.25 edges, whereas the
Corollary~\ref{cor:penalty_bounds}-based T1 and T2 encodings --- despite having
the sparsest $k$-subgraphs as their global minimizers --- collapse to 2.2\% optimality with
mean errors above~70.  The SD baseline reproduces the same encoding
effect without any quantum hardware: AL-SD reaches 92.4\%, T1-SD
9.8\% and T2-SD 8.0\%.  The encoding therefore largely determines the outcome for
the classical baseline just as it does for the QPU.

Meanwhile it appears that with the given settings the QPU
does not provide better starting points than random states.  With
the same 500 reads, AL-SD solves 92.4\% of the instances against
83.1\% for AL-QPU+PP, and under T1 and T2 the mean error after QPU starts (70.75 and
71.06) is about three times that after random starts (24.62 and
24.44).  Raw QPU samples
(AL-QPU, 0.4\%) are far from optimal; nearly all of the solution
quality of QPU+PP is contributed by the descent step.  This differs
from D-Wave's own example, a random Ising problem defined directly on
the QPU's native graph, where QPU starts gave clearly lower energies
than random starts~\cite{DWave_PP_Greedy}.  Our problems are dense and
require chains of up to 32 qubits, which points to the embedding as a
likely cause of the difference.

The penalty table (upper part of Table~\ref{tab:sks_penalties}) shows why the
encoding matters: the AL variants terminate with $\eta<1$ on average,
whereas T1 and T2 present $\eta$ of order 20--50.  With such penalties
every single-flip descent ends at a feasible point (all T1-SD and
T2-SD runs returned cardinality~$k$ in one iteration) but cannot cross
the barriers of Observation~\ref{obs:barriers} between feasible
points, so the result is decided by where the descent starts.
Consistently, T1-SD and T2-SD solve 48.9\% and 40.0\% of the $n=40$
instances and none with $n\ge 80$, while AL-SD solves 84--100\% at
every size.

On the QPU, the penalty is not the only scale in the programmed
Hamiltonian.  Each logical variable is a chain of physical qubits held
together by the chain coupling, which is itself a penalty: it must
dominate the problem terms acting on the chain for the embedded Ising
problem to share its ground state with the logical QUBO
(Section~\ref{sec:distortion_qa}).  A standard sufficient value for
variable~$i$ is $|h_i|+\sum_j|J_{ij}|$~\cite{Choi_2008}, which for a
dense penalised QUBO grows linearly in $n$, in $\eta$ and in
$\lambda$.  Table~\ref{tab:qpu_chains} reports the chain statistics of
the raw AL-QPU runs.  At $n=40$ the fixed chain strength of~15
exceeds this bound and no chain broke; at $n=140$ it is only $0.41$ of
the bound and 95\% of the reads contain at least one broken chain.
Across all runs the fraction of broken reads is predicted by the ratio
of the used to the sufficient chain strength (Spearman $\varrho_{\mathrm{S}}=-0.88$),
and within each size it increases with $\lambda$
($\varrho_{\mathrm{S}}=0.53$--$0.66$).

The chains also set the hardware scale.  Because the chain coupling
value of 15 exceeded every AL problem coupling, the auto-scaling
divided the whole Hamiltonian by~15, and the difference between an
edge and a non-edge was programmed as $0.25/15\approx 1.7\%$ of the
coupler range for every AL run, instead of 11--34\% for the problem
alone.  Raising the chain strength to the sufficient level would
reduce it further, to 0.7--2.5\% for AL and to 0.05\% for T2 at
$n=40$.  The fixed chain strength therefore trades chain integrity
against resolution: the embedding introduces a second penalty whose
magnitude must be set with the same care as $\eta$, and which
interacts with it.

These effects limit what our QPU data can attribute to the hardware.
The 0\% vs.\ 78.7\% contrast between T2-QPU+PP and AL-QPU+PP on the 75
instances with $k=n/2$ is reproduced by the SD baseline (T2-SD 0\%,
AL-SD 90.7\%), so it is mainly an encoding effect, not a precision
mismatch.
Chain breaks alone do not explain the poor raw QPU samples either: at
$n=40$ no chain broke, yet only 1 of 45 raw AL-QPU runs was optimal.
Chain-break statistics were recorded only for the raw runs; the
post-processed runs used the same embedding, the same chain strength
and similar terminal penalties (upper part of Table~\ref{tab:sks_penalties}),
so we take the raw runs as representative of the sampling stage.
Separating the contributions of chain breaks, resolution, noise and
anneal dynamics requires experiments with an $\eta$-dependent chain
strength, which we leave for future work.

Simulated annealing (SA) results are presented in
Table~\ref{tab:sks_sa_aggregate} and the middle part of
Table~\ref{tab:sks_penalties}.  The AL variants (AL-SA) achieve 100\%
optimality across all 270 instances with zero error, confirming that
the augmented Lagrangian loop can reliably find the optimum when
paired with a suitable inner solver. The gap between AL-QPU+PP
(83.1\%) and AL-SA (100\%) therefore reflects the solver, not the
penalty strategy: AL-SD, with the same 500 reads as AL-QPU+PP, lies in
between (92.4\%), and SA explores the landscape far more thoroughly
than a single descent.

The T1 and T2 SA variants, however, perform poorly: T1-SA(10) reaches
only 1.9\% optimality, and even with 1000~reads T1-SA(1000) manages only 8.5\%.
Increasing \texttt{num\_reads} from 10 to 1000 roughly
halves the mean error for both T1 and T2, yet none approaches the AL
loop's 100\%, and T2 is worse than T1 at every \texttt{num\_reads}
level, consistent with $\eta_{T2} = 2\eta_{T1}$.  As SA solves the
QUBO exactly as specified, this degradation is not a hardware effect;
it is consistent with the barriers of Observation~\ref{obs:barriers},
which grow with~$\eta$ and trap the search in suboptimal feasible
basins, as discussed above for SD.  Post-processing the SA samples
with steepest descent changed the optimality rates by less than one
percentage point, since SA samples are already close to single-flip
local minima; we therefore omit the SA+PP variants from
Table~\ref{tab:sks_sa_aggregate}.

The penalty table (Table~\ref{tab:sks_penalties}) confirms
that the Corollary~\ref{cor:penalty_bounds}-prescribed penalties are instance-determined
and solver-independent: they do not depend on \texttt{num\_reads},
and on the 225 shared instances the T1 and T2 values used for SA, SD
and QPU+PP coincide.

Greedy, despite using no penalty formulation at all, achieves
17.8\% optimality with a median error of only~3 ---
outperforming every T1 and T2 variant.
This underscores that the
Corollary~\ref{cor:penalty_bounds}-based penalties, although they
guarantee the correct global minimizers, can be counterproductive when
the inner solver is a heuristic.

\begin{table}[htbp]
    \caption{Aggregate performance of the QPU variants and of the
        steepest-descent baseline (SD, random starts) on the 225
        instances with $n\le 140$, all with 500 reads.  Absolute error
        is $\hat{m}_k - m_k^*$, the number of edges above the optimum
        in the returned $k$-subgraph.  QPU+PP and SD use the same
        descent routine; they differ only in the starting states.}
    \label{tab:sks_qpu_aggregate}
    \centering
    \small
    \begin{tabular}{lrrr}
        \toprule
        Variant   & \% optimal & Mean abs.\ err.\ & Median abs.\ err.\ \\
        \midrule
        AL-QPU    &  0.4 & 50.00 & 45 \\
        AL-QPU+PP & 83.1 &  0.25 &  0 \\
        T1-QPU+PP &  2.2 & 70.75 & 61 \\
        T2-QPU+PP &  2.2 & 71.06 & 64 \\
        \midrule
        AL-SD     & 92.4 &  0.09 &  0 \\
        T1-SD     &  9.8 & 24.62 & 15 \\
        T2-SD     &  8.0 & 24.44 & 14 \\
        \bottomrule
    \end{tabular}
\end{table}

\begin{table}[htbp]
    \caption{Chain statistics of the raw AL-QPU runs (45 runs per size,
        last AL iteration, chain strength 15).  ``Sufficient'' is
        $\max_i(|h_i|+\sum_j|J_{ij}|)$ of the logical Ising problem;
        CBF is the chain-break fraction; ``Edge signal'' is the
        programmed difference between an edge and a non-edge as a
        fraction of the coupler range.}
    \label{tab:qpu_chains}
    \centering
    \small
    \begin{tabular}{rrrrrrr}
        \toprule
        $n$ & Chain len. & Qubits & Reads w/ break & Mean CBF
            & Used/suff. & Edge signal \\
        \midrule
         40 &  4.6 &  185 &  0\% & 0.000 & 1.50 & 1.7\% \\
         80 & 10.1 &  810 &  4\% & 0.001 & 0.71 & 1.7\% \\
        100 & 13.5 & 1349 & 26\% & 0.004 & 0.56 & 1.7\% \\
        120 & 16.3 & 1954 & 70\% & 0.033 & 0.44 & 1.7\% \\
        140 & 20.3 & 2842 & 95\% & 0.084 & 0.41 & 1.7\% \\
        \bottomrule
    \end{tabular}
\end{table}

\begin{table}[htbp]
  \caption{Terminal penalty coefficients. The QPU variants and the SD
    baseline are above the first separator line; the middle part holds
    the SA variants, whereas the bottom part has the MQLib and SBM
    variants. The AL variants grow $\eta$ iteratively until
    cardinality is satisfied. T1 and T2 terminate in a single
    iteration. The T2 encoding is purely quadratic
    ($\lambda \equiv 0$), while T1 operates with equal linear and
    quadratic penalties, so that $\eta_{T2} = 2\eta_{T1}$
    instance-by-instance. (For the minima of $\eta$, because of the
    freedom in the choice of the small value the factor 2 does not
    manifest.) The T1 and T2 penalties are determined by the instance,
    the choice of $\tilde m_k$, and the small positive value of
    $\eta$ used where the boundary value is zero
    (Section~\ref{sec:cs_variants}), so all the sub-variants within each encoding can have
    identical penalty values. Hence, in the middle and bottom parts,
    T1 and T2 rows are collapsed. The GLO variants share the same
    greedy $\lambda_0$ initialization. Their terminal penalties differ
    only through solver-dependent iteration trajectories, and their
    aggregate statistics are nearly identical.}
    \label{tab:sks_penalties}
    \centering
    \small
    \begin{tabular}{lrcrcc}
        \toprule
        Variant & $\lambda$ & $\lambda$
            & $\eta$ & $\eta$ & Mean \\
                & mean & [min, max]
            & mean & [min, max] & iters \\
        \midrule
        \multicolumn{6}{l}{\emph{$n\le 140$ (225 instances)}} \\
        AL-QPU    & 24.0 & $[0,\ 88.1]$  &  0.09 & $[0.05,\ 0.79]$  & 6.03 \\
        AL-QPU+PP & 22.7 & $[0.0,\ 81.0]$  &  0.06 & $[0.05,\ 0.37]$  & 2.32 \\
        T1-QPU+PP & 23.1 & $[0,\ 52.0]$ & 23.1  & $[0.05,\ 52.0]$  & 1.00 \\
        T2-QPU+PP &  0.0 & $[0,\ 0]$       & 46.2  & $[0.05,\ 104.0]$ & 1.00 \\
        AL-SD     & 22.7 & $[0.0,\ 81.0]$  &  0.07 & $[0.05,\ 0.45]$  & 2.68 \\
        T1-SD     & 23.1 & $[0.0,\ 52.0]$  & 23.1  & $[0.05,\ 52.0]$  & 1.00 \\
        T2-SD     &  0.0 & $[0,\ 0]$       & 46.2  & $[0.05,\ 104.0]$ & 1.00 \\
        \midrule
        \multicolumn{6}{l}{\emph{all 270 instances, SA variants}} \\
        AL-SA        &  25.4 & $[0,\ 91.9]$     &  0.12 & $[0.05,\ 3.3]$   & 3.77 \\
        T1-SA        &  25.9 & $[0,\ 59.5]$     &  25.9 & $[0.05,\ 59.5]$  & 1.00 \\
        T2-SA        &     0 & $[0,\ 0]$        &  51.7 & $[0.05,\ 119.0]$ & 1.00 \\
      \midrule
      \multicolumn{6}{l}{\emph{all 270 instances, MQLib and SBM}} \\
        T1 (HH, SBM) &  25.9 & $[0,\ 59.5]$  & 25.9 & $[0.50,\ 59.5]$  & 1.00 \\
        T2 (HH, SBM) &     0 & $[0,\ 0]$     & 51.7 & $[0.50,\ 119.0]$ & 1.00 \\
        GLO/MQ       &  25.4 & $[0,\ 91.9]$  &  0.14 & $[0.10,\ 0.81]$ & 2.47 \\
        GLO/SBM(500) &  25.4 & $[0,\ 91.9]$  &  0.13 & $[0.10,\ 0.81]$ & 2.02 \\
        \bottomrule
    \end{tabular}
\end{table}
\begin{table}[htbp]
    \caption{Aggregate performance of SA and Greedy.
        (270 instances).}
    \label{tab:sks_sa_aggregate}
    \centering
    \small
    \begin{tabular}{lrrr}
        \toprule
        Variant   & \% optimal & Mean abs.\ err.\ & Median abs.\ err.\ \\
        \midrule
        AL-SA            & 100.0 &  0.00 & 0    \\
        T1-SA(10)        &   1.9 & 58.12 & 49.5 \\
        T1-SA(1000)      &   8.5 & 39.73 & 32   \\
        T2-SA(10)        &   1.9 & 75.37 & 68   \\
        T2-SA(1000)      &   4.4 & 54.97 & 47.5 \\
        Greedy           &  17.8 &  3.97 & 3    \\
        \bottomrule
    \end{tabular}
\end{table}

\subsection{MQLib and SBM}
\label{sec:cs_mqlib}

In this subsection we evaluate the penalty strategies using two inner
solvers that are both capable of handling large and dense QUBO
instances, and do not require embedding. MQLib's hyper-heuristic (HH),
as a portfolio of classical QUBO heuristics, is typically robust and
efficient - it finds and uses the best-fitting traditional classical
heuristic. We use MQLib setting a fixed wall-time budget (1\,s or
10\,s per iteration, depending on the experiment), and HH has no other
hyperparameters to tune.

As opposed to that, simulated bifurcation (SBM), despite being
classical, works by a different mechanism than the heuristics covered
by MQLib: a continuous-variable dynamical system evolves real-valued
amplitudes under a smooth ordinary differential equation, and binary values are obtained (via the
sign) only at readout. We use one fixed setting of hyperparameters in
all experiments: the
continuous-time variant with 500 reads, 2000 steps and
$\Delta t = 0.1$.  This setting was chosen ad hoc and not tuned.
Because every penalty strategy is run with the same setting, the
comparison between strategies remains within-solver; however, the
absolute performance of SBM, and possibly the size of the gap between
strategies, may change with tuning.  In particular, a large~$\eta$
changes the scale and spectrum of the QUBO, so hyperparameters suited
to the small penalties of the AL loop need not suit the
Corollary~\ref{cor:penalty_bounds}-prescribed ones.  As with the chain
strength on the QPU (Section~\ref{sec:cs_results}), the solver's own
parameters thus form a further layer at which the encoding interacts
with performance.

A large $\eta$ is likely detrimental to SBM for related reasons. SBM's
dynamics evolve continuous oscillator amplitudes under a bifurcation
process governed by the dominant mode of the Ising coupling matrix --
the same $J_{ij}, h_i$ obtained from the QUBO-to-Ising conversion
described in Section~\ref{sec:distortion_qa}. The quadratic penalty
contributes a uniform, rank-one term to this matrix, which grows with
$\eta$ and is independent of $Q$. Once this contribution comes to
dominate the matrix, it is plausible that the dynamics are steered by
the penalty structure rather than by the underlying graph, in a way
related to how the hardware-precision mechanism of
Section~\ref{sec:distortion_qa} lets the same uniform term dominate
the objective's contribution to the QPU couplings. We have not
investigated this mechanism in detail, and leave a proper
characterization to future work.

We have evaluated all the 270 benchmark instances.  As in
Section~\ref{sec:cs_results}, the T1 and T2 penalties use the known
optimum~$m_k^*$.  Concerning performance, let us analyze the data in
Table~\ref{tab:sks_hh_aggregate}.  The pattern is consistent with the
QPU/SA findings: the GLO variants (GLO/MQ and GLO/SBM) achieve 100\%
optimality, while the Corollary~\ref{cor:penalty_bounds}-based T1 and
T2 encodings show clear degradation.  The MQLib hyper-heuristic is
robust enough to largely compensate for the distortion on this
benchmark --- T1/HH still reaches 97.4\% --- but the gap is visible
and widens for T2/HH (94.4\%).  The SBM solver with the given fixed
settings shows significantly worse performance under the
Corollary-based penalties: T1/SBM converges on only 41 of 270
instances (15.2\%) and T2/SBM on 32 (11.9\%).

Penalties are tabulated in the lowest part of Table~\ref{tab:sks_penalties}.
The convergence failure is not a solver defect but a direct
consequence of the encoding: T1/SBM and T1/HH receive identical
penalty coefficients (both instance-determined from the
Corollary~\ref{cor:penalty_bounds} with $m_k^*$), yet MQLib finds
near-optimal solutions under the same $\eta \approx 26$ that the SBM
cannot handle with its current settings. This illustrates that the
distortion effect is universal in direction --- all heuristic solvers
degrade with increasing~$\eta$ --- but its severity depends on the
solver's robustness against the negative consequences of big $\eta$.
The penalty table also reveals a sharp quantitative contrast
between the two strategy families.  The Corollary~\ref{cor:penalty_bounds} prescribes
$\eta \approx 26$ for T1 and $\eta \approx 52$ for T2; the GLO
variants terminate with $\eta < 1$ on average (max $0.81$), a
factor of 30--60$\times$ smaller.  

The maximal number of iterations differ between the variants, which
deserves comment.  The T1/SBM and T2/SBM variants are limited to
\texttt{max\_iter}~=~1, meaning the algorithm cannot retry if the
first inner solve fails to produce a feasible cardinality.  This is by
design: the sparsest $k$-subgraphs are global QUBO minimizers already
on the first call (Corollary~\ref{cor:sks_boundary}), so additional iterations would only
help if the solver returned a better approximate solution --- which,
given values as high as $\eta \approx 26$--$52$, is unlikely to
be the case.  The GLO variants, by contrast, are allowed up to 100
iterations and use the AL loop to steer the solver toward feasibility
through small, iterative adjustments.  The mean iteration counts of
2.0--2.5 confirm that the loop converges quickly once $\lambda$ is
near its optimal value.

\begin{table}[htbp]
    \caption{Aggregate performance of the six algorithm variants
        over the 270 benchmark instances.  Absolute error is
        $\hat{m}_k - m_k^*$.  Non-converged instances are excluded
        from the error and optimality statistics; T1/SBM and
        T2/SBM converged on only 41 and 32 of 270 instances,
        respectively.}
    \label{tab:sks_hh_aggregate}
    \centering
    \small
    \begin{tabular}{lrrrr}
        \toprule
        Variant & Conv.\ & \% optimal
            & Mean abs.\ err.\ & Median abs.\ err.\ \\
        \midrule
        T1/HH        &     270 &  97.4 &  0.09 & 0   \\
        T1/SBM       &  41/270 &  41.5 & 11.93 & 3   \\
        T2/HH        &     270 &  94.4 &  2.29 & 0   \\
        T2/SBM       &  32/270 &  28.1 &  8.59 & 7.5 \\
        GLO/MQ       &     270 & 100.0 &  0.00 & 0   \\
        GLO/SBM(500) &     270 & 100.0 &  0.00 & 0   \\
        \bottomrule
    \end{tabular}
\end{table}

\section{Penalty Bounds: Comparison with Lasserre and EXPEDIS}
\label{sec:cs_cor}

Let us now compare the quadratic penalty prescribed by
Corollary~\ref{cor:penalty_bounds} with two general-purpose bounds
from the literature, Lasserre's approach~\cite{Lasserre2016} and
EXPEDIS (Gusmeroli and Wiegele~\cite{EXPEDIS}).  Both are
pure-quadratic approaches ($\lambda = 0$), hence, we first compare
against our T2 bound ($\lambda = 0$), and then quantify the additional
improvement from T1 ($\lambda = \eta$), which shifts a part of the
penalisation to the distortion-free linear term.

Lasserre~\cite{Lasserre2016} first rewrites the problem in the
variables $\tilde{x} = 2x - \mathbf{e} \in \{-1,1\}^n$, in which the
objective reads $\tilde{c}^\top \tilde{x} + \tilde{x}^\top \tilde{F}
\tilde{x}$ up to an additive constant; for the SkS problem
$\tilde{F} = A/8$ and $\tilde{c} = A\mathbf{e}/4$.  His penalty
coefficient is $\eta_{\mathrm{Las}} = 2\rho_{\mathrm{Las}} + 1$, where
$\rho_{\mathrm{Las}} := \max\{|r^1|,\,|r^2|\}$, and $r^1$ and $r^2$
are the minimum and the maximum of the standard (Shor) semidefinite programming (SDP) relaxation
of this objective over $\{-1,1\}^n$, with the linear constraints
ignored.  Written in the original variables, the penalty term
of~\cite{Lasserre2016} is exactly $(\mathbf{e}^\top x - k)^2$, so
$\eta_{\mathrm{Las}}$ is directly comparable with our~$\eta$.  In
general, computing $\rho_{\mathrm{Las}}$ requires solving two SDPs.
For the SkS problem all entries of $\tilde{c}$ and $\tilde{F}$ are
nonnegative, and every feasible point $(x, X)$ of the relaxation has
$|x_i| \leq 1$ and $|X_{ij}| \leq 1$; hence $|r^1| \leq r^2 =
\mathbf{e}^\top \tilde{c} + \mathbf{e}^\top \tilde{F} \mathbf{e} =
\tfrac{3}{4}|E|$, where $|E|$ is the number of edges of~$G$.  This
gives $\eta_{\mathrm{Las}} = \tfrac{3}{2}|E| + 1$, which grows with
graph size and density~$\bar{\rho}$ regardless of~$k$.

EXPEDIS tightens Lasserre's bound by separating feasible from
infeasible objective values:
$\eta_{\mathrm{GW}} = u - \ell + \epsilon$ (denoted $\sigma$
in~\cite{EXPEDIS}), where $u$ is an upper bound on
$u^* = \max_{\mathbf{e}^\top x = k} f(x)$, the densest $k$-subgraph
value of~$G$, and $\ell$ is a lower bound on
$\ell^* = \min_{\mathbf{e}^\top x \neq k} f(x)$; both bounds are
obtained from SDP relaxations.  For our problem $f(x) \geq 0$ for every $x$, and $f(\mathbf{0}) = 0$
at the infeasible point $x = \mathbf{0}$; hence $\ell^* = 0$.  We
therefore take $\ell = 0$, the best possible lower bound, and
$\eta_{\mathrm{GW}} = u_\Delta + \epsilon$, where $u_\Delta \geq u^*$
is the SDP bound of~\cite{EXPEDIS} that includes the cardinality
constraint.

As opposed to these, our T2 bound exploits the cardinality structure
directly: $\eta_{T2} = \min\{\tilde{m}_k,\,k{-}1\} + \epsilon$, with
$\lambda = 0$ (the experiments of
Section~\ref{sec:computational_results} use the boundary value
$\epsilon = 0$, see Section~\ref{sec:cs_variants}).  The $\min$ operator is the key structural advantage:
for dense graphs or large~$k$, the bound saturates at $k{-}1$ rather
than growing with the edge count.  It can be evaluated with a greedy
pass, no SDP is needed.

From the definitions one obtains a hierarchy of the penalty values.
Disregarding the arbitrarily small~$\epsilon$, the three
pure-quadratic bounds satisfy
\begin{equation}\label{eq:penalty_hierarchy}
  \eta_{T2}
  \;\leq\; \eta_{\mathrm{GW}}
  \;<\; \frac{2}{3}\,\eta_{\mathrm{Las}}.
\end{equation}
The left inequality follows from
$\eta_{T2} = \min\{\tilde{m}_k,\,k{-}1\} \leq \tilde{m}_k \leq u^* \leq
u_\Delta = \eta_{\mathrm{GW}}$, since no $k$-subgraph has more edges
than the densest one; it therefore holds for every feasible
estimate~$\tilde{m}_k$, in particular for $m_k^*$ and for the greedy
value.  The right inequality follows from
$u_\Delta \leq |E| = \tfrac{2}{3}(\eta_{\mathrm{Las}} - 1)$: the
objective of the relaxation defining $u_\Delta$ cannot exceed the sum
of its nonnegative coefficients, which is~$|E|$.  On all our benchmark
instances even $\eta_{\mathrm{GW}} \leq \eta_{\mathrm{Las}}/2$ holds.

The gap between T2 and EXPEDIS is not a constant factor: when
$k{-}1 \ll m_k^*$, the $\min$ in the T2 bound saturates at
$k{-}1$ while $\eta_{\mathrm{GW}}$ is at least~$m_k^*$,
so the ratio grows without bound.  In the extreme case where the
optimal $k$-subgraph is an independent set ($m_k^* = 0$), our
framework requires only $\eta > 0$ while Lasserre prescribes
$\eta_{\mathrm{Las}} = \tfrac{3}{2}|E| + 1$: the improvement is qualitative,
not merely quantitative.

Our T1 bound provides further improvement by allowing for the linear penalty, setting
$\lambda_{T1} = \eta_{T1} = \min\{\tilde{m}_k/2,\,(k{-}1)/2\}$.  The full hierarchy is
therefore
\begin{equation}\label{eq:full_hierarchy}
  \eta_{T1}
    \;=\; \frac{\eta_{T2}}{2}
  \;\leq\; \frac{\eta_{\mathrm{GW}}}{2}
  \;<\; \frac{\eta_{\mathrm{Las}}}{3}.
\end{equation}
As
established in Sections~\ref{sec:landscape_roughening}
and~\ref{sec:distortion_qa}, the linear penalty does not contribute to
landscape roughening or hardware precision loss. 
Table~\ref{tab:penalty_comparison_summary} summarizes the four
approaches.  Each row tightens the bound further, trading
generality for problem structure.
\begin{table}[htbp]
  \caption{Hierarchy of quadratic penalty bounds for the SkS
    problem.  Each successive approach exploits more problem
    structure, yielding a tighter~$\eta$ at comparable or lower
    computational cost.}
  \label{tab:penalty_comparison_summary}
  \small
  \centering
  \begin{tabular}{llll}
    \toprule
    Approach & $\eta$ (or $\sigma$) & Computation
      & $\lambda$ \\
    \midrule
    Lasserre~\cite{Lasserre2016}
            & $\tfrac{3}{2}|E| + 1$
      & two SDPs in general; closed form for SkS
      & $0$ \\[3pt]
    EXPEDIS~\cite{EXPEDIS}
            & $u_\Delta + \epsilon$
      & SDP
      & $0$ \\[3pt]
    Cor.~\ref{cor:penalty_bounds} (T2)
      & $\min\{\tilde{m}_k,\,k{-}1\} + \epsilon$
      & greedy pass
      & $0$ \\[3pt]
    Cor.~\ref{cor:penalty_bounds} (T1)
      & $\frac{\min\{\tilde{m}_k,\,k{-}1\}}{2} + \epsilon$
      & greedy pass
      & $= \eta_{T1}$ \\
    \bottomrule
  \end{tabular}
\end{table}

We evaluated the three pure-quadratic bounds on the 270 benchmark
instances, with $\tilde{m}_k$ taken from the greedy pass and
$\eta_{\mathrm{GW}} = \lfloor u_\Delta \rfloor + 1$.  On the 265
instances with $\eta_{T2} > 0$, the T2 bound is smaller than the
EXPEDIS penalty by a factor of 4--47 (median~13) and smaller than
Lasserre's by a factor of 10--367 (median~90); both factors grow with
$n$ and with the density~$\bar{\rho}$, and they double for T1.  This
demonstrates the benefit of exploiting the problem structure when
calculating penalty coefficients.

\section{Summary and Outlook}
\label{sec:cs_conclusion}

We have studied the encoding of quadratic binary optimization problems
with linear constraints into QUBO form, on the example of the sparsest
$k$-subgraph problem, which has a single cardinality constraint. Our results underline the widely
accepted view that the inappropriate application of the commonly
used quadratic penalty approach can have significant negative effects
on the solution quality. We consider in detail two inner-solver-dependent
mechanisms (applicable to simulated annealing and quantum
annealing, respectively) that explain why large quadratic penalties
are harmful, while our extensive experiments support the conjecture that
this effect arises with all kinds of QUBO solvers and probably
for more general constraints, too.

For the particular case studied, we identify three ways of mitigating
the problems caused by quadratic penalties. We derive sufficient bounds
for the penalty magnitudes based on the structure of the particular
constraint, which lead to significantly smaller penalties than those
resulting from general-purpose approaches. We argue that if the problem's
structure allows for introducing a linear penalty term in addition to
the quadratic one, as in the problem studied here, this can be used
to decrease the quadratic penalty while not contributing to landscape
distortion or hardware precision issues. Finally, our results demonstrate
that iterative determination of the penalty coefficients, by repeated
application of the inner solver, can be the most successful strategy.

Our work is not limited to a single approach to solving QUBOs: we have
used steepest descent from random starts and
simulated annealing as classical heuristics, MQLib's classical
hyper-heuristic, a D-Wave quantum annealer, as well as the simulated
bifurcation machine (SBM), a recently introduced quantum-inspired
solver. This holistic
approach to inner QUBO solvers confirms that in all cases, the
encoding of the constraints is as important as the choice or
performance of the solver. The experiments are instructive
from the point of view of the behavior of these solvers. Meanwhile,
the results suggest that some of the phenomena we have observed are
not specific to the problem or to the inner solver. The generalization to a broader
class of constraints and inner solvers, and the identification of
possible deeper principles behind these observations, are left to future
research.

\begin{acknowledgements}
  The authors gratefully acknowledge the J\"ulich Supercomputing
  Centre (\url{https://www.fz-juelich.de/jsc}) for funding this project
  by providing computing time on the D-Wave Advantage\texttrademark{}
  System JUPSI through the J\"ulich UNified Infrastructure for Quantum
  computing (JUNIQ).
  
\end{acknowledgements}

\section*{Declarations}

\paragraph{Funding}
This research was co-funded by the Slovenian Research and Innovation Agency ARIS through the research project Quantum Solver for Hard Binary Quadratic Problems (Grant No. J7-50186) and through the annual work programme of Rudolfovo; and by the
QEC4QEA project (Grant Agreement No. 101194322), funded by the European High Performance Computing Joint Undertaking (EuroHPC JU) and participating countries, including Slovenia, under the Horizon Europe research and innovation programme.

\paragraph{Competing interests}
The authors have no relevant financial or non-financial interests to
disclose.

\paragraph{Data availability}
The datasets generated and analyzed during the current study are
available from the corresponding author on reasonable request.

\paragraph{Code availability}
Code and data are available from the corresponding author on
reasonable request.

\paragraph{Author contributions}
The original ideas are mostly due to R.~Ku\v{z}el, with contributions
from J.~Povh. R.~Ku\v{z}el carried out the analytical derivations, and
M.~Koniorczyk refined and verified the theoretical results.
R.~Ku\v{z}el and M.~Koniorczyk jointly designed and performed the
numerical experiments and developed the interpretation of the results.
M.~Koniorczyk wrote the manuscript with input from R.~Ku\v{z}el.
J.~Povh supervised the project and proofread the manuscript. All
authors read and approved the final manuscript.

\paragraph{Use of generative AI}
During the preparation of this work the authors used Claude
(Anthropic) for language editing, literature search, coding, and
analysis.  The research ideas belong to the authors, and Claude was
used as an assistant under human supervision.  The authors reviewed
and edited all output and take full responsibility for the content of
the publication.


\end{document}